\documentclass[preprint,review,12pt]{elsarticle}
\usepackage{amssymb}
\usepackage{lineno}
\usepackage{graphicx}
\usepackage{tabularx}
\usepackage{url}
\usepackage[american]{babel}
\usepackage{amsmath,amsfonts}
\usepackage{algorithmic}
\usepackage{graphicx}
\usepackage{textcomp}
\usepackage{xcolor}
\usepackage{mathtools}
\usepackage{booktabs} 
\usepackage{tikz} 
\usepackage{footmisc}
\usepackage{subcaption}
\usepackage{booktabs}
\usepackage{siunitx}
\usepackage[makeroom]{cancel}
\usepackage{MnSymbol}
\usepackage{algorithm}
\usepackage{float}
\usepackage{dblfloatfix}
\usepackage[colorlinks]{hyperref}

 \newcommand{\argmin}{\operatornamewithlimits{argmin}} 
\newcommand{\argmax}{\operatornamewithlimits{argmax}}

\newcommand{\Pro}{{I\kern-.35em P}}
\newcommand{\tinysquare}{\scriptscriptstyle{\square}}

\DeclarePairedDelimiter\ceil{\lceil}{\rceil}

\newtheorem{theorem}{Theorem}

\newtheorem{definition}{Definition}

\journal{Reliability Engineering and System Safety\, }

\begin{document}

\begin{frontmatter}



\title{A Data-driven Approach to Risk-aware Robust Design}


\author[inst1]{Luis G. Crespo\corref{cor1}}
\ead{Luis.G.Crespo@nasa.gov}
\cortext[cor1]{Corresponding author.}
\affiliation[inst1]{organization={NASA Langley Research Center}, 
            city={Hampton},
            state={VA},
            postcode={23681}, 
            country={USA.}}
\author[inst1]{Bret K. Stanford}
\author[inst1]{Natalia Alexandrov.}

\begin{abstract}
This paper proposes risk-averse and risk-agnostic formulations to robust design in which solutions that satisfy the system requirements for a set of scenarios are pursued. These scenarios, which correspond to realizations of uncertain parameters or varying operating conditions, can be obtained either experimentally or synthetically. The proposed designs are made robust to variations in the training data by considering perturbed scenarios. This practice allows accounting for error and uncertainty in the measurements, thereby preventing data overfitting. Furthermore, we use relaxation to trade-off a lower optimal objective value against lesser robustness to uncertainty. This is attained by eliminating a given number of optimally chosen outliers from the dataset, and by allowing the perturbed scenarios to violate the requirements with an acceptably small probability.  For instance, we can seek a design that satisfies the requirements for as many perturbed scenarios as possible, or pursue a riskier design that attains a lower objective value in exchange for a few scenarios violating the requirements. These ideas are illustrated by considering the design of an aeroelastic wing.
\end{abstract}


\begin{highlights}
\item Data-driven formulations to robust design are proposed. 
\item The resulting designs are robust to uncertainty and error in the data.
\item Means to trade off performance against robustness are developed. 
\item Strategies to lower the computational cost are proposed and tested.
\item The design of an aeroelastic wing is used for illustration.
\end{highlights}


\end{frontmatter}


\section{Introduction} 
Optimization under uncertainty addresses the challenge of making optimal decisions about a system whose performance depends on parameters having an unknown value. Uncertainty might be caused by aleatory variability, ignorance, or changing operating conditions. Traditional deterministic optimization assumes that all the parameters of such a model are known, a rare case in real-world applications. Optimization under uncertainty, on the other hand, is a decision-making process that accounts for the effects of uncertainty upfront.  A plethora of applications in finance, science and engineering, -including  regression analysis, model calibration, controls, system identification, machine learning, systems theory, financial mathematics, structural design-  require decision-making under uncertainty. 

The mathematical framework supporting the proposed strategies is introduced next. Let $J:\Theta \rightarrow \mathbb{R}$ be an objective function of the decision variable $\theta\in{\mathbb R}^{n_\theta}$.  Furthermore, let $r_k(\theta,\delta): \mathbb{R}^{n_{\theta}} \times \mathbb{R}^{n_{\delta}} \rightarrow  \mathbb{R}$ for $k=1,\ldots n_r$ be a continuous function associated with a system requirement, where the parameter $\delta\in \Delta \subset \mathbb{R}^{n_\delta}$ is subject to aleatory uncertainty. The optimization program under consideration is
\begin{align}
	 \min_{\theta \in \Theta} & \quad J\left(\theta \right) \label{sta}\\
	\text{ subject to:} & \quad r_k(\theta,\delta)\leq 0,  \text{ for } k=1,\ldots n_r \text{ and all }\delta\in\hat{\Delta},\nonumber 
\end{align}
where $\hat{\Delta}\subseteq\Delta$ is yet to be prescribed. The objective and constraints of (\ref{sta}) are assumed to be continuous functions of $\theta$, thereby making standard gradient based algorithms applicable to the forthcoming formulations. 

Denote as $\theta^\star$ the solution to (\ref{sta}). The \emph{feasible set} of (\ref{sta}) corresponding to a parameter point $\delta$ is 
\begin{equation}
\label{feaset}
\mathcal{X}(\delta)=\left \{\theta\in \Theta: r_k(\theta,\delta)\leq 0,\; k=1,\ldots n_r \right\}.
\end{equation}
Hence, $\mathcal{X}(\delta)$ is comprised of all the decision (or design) points that satisfy the constraints for the parameter point $\delta$. The \emph{success domain} of (\ref{sta}) corresponding to the decision point $\theta$ is
\begin{equation}
\mathcal{S}(\theta)=\bigcap_{k=1}^{n_r} \mathcal{S}_k(\theta),\label{su}
\end{equation}
where $\mathcal{S}_k(\theta)=\left \{\delta\in \Delta:\; r_k(\theta,\delta)\leq 0\right\}$. Hence,  $\mathcal{S}(\theta)$ is comprised of all the parameter points for which the decision point $\theta$ satisfies the constraints. The set $\Delta$ is partitioned into the \emph{success} and the \emph{failure} domains, denoted as $\mathcal{S}(\theta)$ and $\mathcal{F}(\theta)$ respectively. These sets satisfy $\mathcal{F}(\theta)\cup \mathcal{S}(\theta)=\Delta$ and $\mathcal{F}(\theta)\cap \mathcal{S}(\theta)=\emptyset$. Each of the sets in the right-hand side of (\ref{su}) is an individual success domain, with their complement being the corresponding individual failure domain.  

The field of optimization under uncertainty was pioneered in the 1950s by Dantzig \cite{Dantzig55} and Charnes \cite{Charnes59}, who set the foundations for \emph{Robust Optimization} (RO) and \emph{Stochastic Optimization} (SO), respectively. 

RO \cite{Soyster73,Kall10,Bertsimas14} solves (\ref{sta}) for a given $\hat{\Delta}$. The main advantage of such a design is the guarantee that the requirements are satisfied for all the realizations of $\delta$ in $\hat{\Delta}$.  However, RO programs are generally intractable even when the requirement functions are convex, e.g., the robust counterpart of a second order cone program with polyhedral uncertainty is NP hard\footnote{Tractable linear RO problems, which require a suitably chosen $\Delta$ and particular forms for $r_k(\theta,\delta)$, are summarized in \cite{Ben-Tal09, Bertsimas11}.} \cite{Elghaoui98,Ben-Tal00,Bertsimas04}.
RO programs are driven by the worst-case combination of uncertainties in $\hat{\Delta}$. As such, other designs can outperform RO designs for most elements of $\hat{\Delta}$.  

SO prescribes the uncertainty in $\delta$ by the probability distribution $\mathbb{P}_\delta$ with support set $\Delta$. This distribution is inferred from data, expert opinion or a mixture of these. Distributions enable the analyst to consider moments of a response function \cite{Tsatsanis98,Mehlawat21,Hammond24} as well as \emph{Chance-Constraints} (CC). SO programs subject to CCs are solved by \emph{Chance-Constrained Optimization} (CCO) \cite{Shapiro09,Coit19}. In contrast to RO, CCO replaces the constraints in (\ref{sta}) with $\mathbb{P}_\delta[{\mathcal{S}_k(\theta)}]\geq \gamma_k$ for $0\ll\gamma_k< 1$. This CC makes the system satisfies the requirement for most uncertainty realizations of $\delta$ on $\Delta$ without prescribing $\hat{\Delta}$ explicitly. CCO designs reduce the intrinsic conservatism of RO designs by ignoring the worst-case realizations of $\delta$ in $\Delta$, $\Delta\setminus \hat{\Delta}$, thereby lowering the objective function. The main difficulty in solving CCO programs stems from the inability to evaluate CCs accurately.  This task entails calculating a high-dimensional integral over a complex integration domain, a problem that is known to be NP-hard \cite{Luedtke10, Qiu14}.  Furthermore, the feasible set of a CCO program is often non-convex, making the identification of $\theta^\star$ difficult.  In fact, the feasible set of a CC corresponding to a convex constraint, i.e., $r_k(\theta, \delta)$ is a convex function of $\theta$ for any $\delta$ and all $k=1,\ldots n_k$, might be non-convex\footnote{The convexity of the requirement functions is preserved by reformulating the CC in terms of the {\em Conditional Value at Risk} (CVaR) at the expense of reducing the feasible set \cite{Shapiro09,Chapman22}. In contrast to CCs, however, CVaR-based constraints make the resulting design risk-averse, e.g., $\theta^\star$ depends on a risk measure taking values on the set $\{\delta: r(\theta^\star,\delta)>0\}$.  
}. CCO programs taking forms that make the feasible set convex can be efficiently solved \cite{Lagoa99, Calafiore06, Henrion08, Henrion11, Prekopa11, Van15}. However, many problems in science and engineering applications are subject to implicit and non-convex requirement functions requiring numerical simulation for their evaluation. In spite of the above drawbacks, the ability to relax the robustness specifications enforced by RO, thereby to reduce the potentially high conservatism of a design synthesized for the worst-case combination of uncertainties, makes CCO an important and practical tool for decision making under uncertainty \cite{Geng19,Mammarella22}. 

A particular case of CCO, called \emph{Reliability-Based Design Optimization} (RBDO), seeks designs that either minimize or bound the probability of failure  \cite{COIT2019106259}. 
A standard approach to RBDO involves two nested loops: an outer loop that searches for the optimal design, and an inner loop that evaluates $\mathbb{P}_\delta[{\mathcal F}(\theta)]$ for every candidate design chosen by the outer loop \cite{ENEVOLDSEN1994169}. This evaluation is often called a reliability analysis.  The high computational cost of accurately estimating small probabilities makes RBDO expensive. Decoupling approaches \cite{YUAN2014107}, single-loop methods \cite{SHAN20081218, MENG201995, WANG2020113436}, and efficient approximations to the failure probability have been used to make RBDO more efficient. Single-loop methods combine the outer and inner loops by substituting the reliability analysis with an approximation \cite{YAO201328}, whereas decoupling methods replace the nested optimization with a sequence of deterministic programs \cite{TORII2019106499}. Some of the approaches used to reduce the computational cost of the inner loop include subset simulation \cite{LI2010384}, line sampling \cite{DEANGELIS2015170}, importance sampling \cite{CHAUDHURI2020106853}, first-order and second-order reliability methods \cite{ENEVOLDSEN1994169, NIKOLAIDIS1988781, Cizelj1994, SCHUELLER2004463}, and many others \cite{DIGE2018431, LI2019106432, Ullmann2015MultilevelEO}. Strategies aiming to reduce the high computational effort of RBDO by using surrogates are available \cite{LI20108966, CHEN2013233, PEHERSTORFER201761,doi:10.1137/17M1122992}. Among the plethora of studies in the literature, we highlight \cite{Zhang21,Ma23}, where a double-loop method and a Gaussian Process Model (GPM) surrogate are combined to speed up computations. In \cite{Van23}, a learning particle swarm optimization approach is decoupled from a GPM surrogate to perform RBDO by sequentially searching for the optimal most probable point. In \cite{Hong24} sequential sampling, Bayesian optimization, and Bayesian reliability analysis are combined to update a GPM surrogate. In \cite{Kera25} an adaptive response surface method is paired with directional sampling to carry out multi-disciplinary RBDO of a launch vehicle. Moreover,  \cite{Chen25} proposes a learning function and point selection strategy to update polynomial chaos-kriging quantile surrogates for RBDO. Finally, \cite{Zhang26} proposes an adaptive conjugate gradient framework and a first-order inverse reliability method to perform RBDO with dependent variables.

An important drawback of the RO and CCO approaches is the high sensitivity of the resulting design to the assumed \emph{uncertainty model}, i.e. $\hat{\Delta}$ and $\mathbb{P}_\delta$ respectively. In the CCO case, the prescription of a distribution generally involves learning a random vector, a dependency structure, and tail model from data. This process is challenging when the number of uncertain parameters, $n_\delta$, is large, when parameter dependencies are strongly nonlinear, when the dataset is incomplete, e.g., the dependency between disjoint subsets of data cannot be measured, or when the data are scarce. Poorly chosen uncertainty models might lead to designs that grossly underperform in practice \cite{ROCKAFELLAR2010499, Sarykalin2008}. Non-probabilistic models and mixtures of non-probabilistic and probabilistic models have been used when the data are scarce \cite{BENHAIM1994227, MENG2020773, CRESPO2019104560, ROCCHETTA2018710}. These models make use of evidence theory \cite{SHAFER20167, ferson2002constructing}, possibility theory \cite{Dubois2001}, credal sets, fuzzy sets and ambiguity sets theory \cite{walley1991statistical, ZADEH1965338, 8693861}.  

The subjectivity in prescribing a distribution can be mitigated by using Distributionally Robust Optimization (DRO) \cite{Rahimian19, Zhengqing21,Lasserre18, Lam16,Bertsimas18, Hanasusanto15,Chen18}. This subjectivity is caused by learning a distribution from limited and possibly corrupted data. DRO frames decision-making as a minimax problem: the decision-maker chooses a policy to minimize a performance criterion, e.g., expected loss or a risk measure such as VaR or CVaR, while an adversary selects a distribution from an ambiguity set to make that criterion as large as possible.  As in RO, this hedges distributional misspecification via a worst-case perspective. DRO is viable when the functional form of the response/requirement functions (or their relaxation) and the ambiguity set yield a tractable maximization in the inner loop. The computational requirements associated with this need depend on the manner by which these functions depend on both $\theta$ and $\delta$, thereby limiting its applicability. DRO has emerged as an important paradigm for machine learning, statistics, controls, finance, engineering, and operations research. 

The ambiguity set is typically categorized as moment-based or discrepancy-based.  In discrepancy-based DRO, ambiguity sets are typically specified as neighborhoods of a reference distribution. From this category we highlight the engineering applications in \cite{Dapogny23,Chen24}, where DRO is used for shape and topology optimization using the Wasserstein distance and the $\phi$-divergence. Furthermore, \cite{Zhang24} uses a distributionally robust CVaR optimization for distribution system planning with short- and long-term uncertainties resulting from extreme weather events, whereas \cite{Ren24} uses DRO to design power-system schedules to recover from weather-caused disasters given an uncertain power demand.

DRO faces several limitations, primarily related to computational complexity, high-dimensional uncertainties, and the difficulty in choosing the ambiguity set and calibrating its parameters \cite{ELDRED20111092, Xie_DRCCP_Wassestein2018, NANNAPANENI20169}. Note that small data sets require a large ambiguity set, which in turn might lead to an intrinsically conservative design that would likely underperform other designs in practice (as it is the case with RO designs). This outcome is a consequence of the worst-case element of the ambiguity set being significantly different from the true distribution. Furthermore, the effective application of DRO to problems combining moments and CCs of arbitrary and possibly implicit requirement functions, such as the ones studied herein, might be intractable or overly-conservative. 
In contrast to RDO, the design framework herein accounts for the effects of measurement noise, limited observability (Footnote 5) and systematic errors in the metrology system, e.g., the uncertainty in the data is prescribed according to the precision and accuracy of the sensor, by modeling and propagating the uncertainty in each datum.

\section{Solution Strategy and Goals}
The strategies proposed below are data-driven. Data-driven methods include the scenario approach \cite{Calafiore06,Campi21,Paulson21,Campi21ml,Garatti22}, the sample average approximation \cite{Luedtke08}, and robust optimization over safe approximations of the feasible CC set \cite{Nemirovski07}.  Developments in the latter two categories, however, have mostly focused on convex requirement functions taking particular forms. The strategies proposed below, which belong to the former category, do not suffer from such restrictions, thereby being broadly applicable. 

Scenario-based strategies make direct\footnote{When such realizations are unavailable, samples drawn from a synthetic distribution could be used instead.} use of a finite number of realizations of $\delta$. This practice eliminates the subjectivity introduced by having to prescribe an uncertainty model from data. More importantly, these strategies offer a computationally cheaper alternative to RBDO approaches while yielding a feasibility guarantee \cite{Garatti22}. These realizations, called the \emph{nominal} scenarios, comprise the dataset 
\begin{equation}
{\mathcal D}=\left \{\delta^{(i)}\right \}_{i=1}^n.\label{ds}
\end{equation}
The nominal scenarios are drawn from an unknown distribution $\mathbb{P}_\delta$.  Hence, we will solve (\ref{sta}) for $\hat{\Delta}={\mathcal D}$ hoping that the resulting design $\theta^\star$ will satisfy the requirements for neighboring scenarios beyond those used for training.

The main goals of this article are two-fold. First, we want to safeguard the robust design against data overfitting. This is of particular interest when the metrology system is inaccurate, e.g., systems commonly used in biological\footnote{Of particular interest are the uncertainties associated with radiation quality effects, as they remain the dominant source of uncertainty for cancer risk projections of astronauts participating in deep space missions.  In this application, the predicted variable of interest, which quantifies the increased biological effectiveness of energetic ions compared to gamma rays, cannot be directly measured. Instead, raw biological measurements, such as tumor incidence and chromosome aberrations, are post-processed with well-documented methods in order to quantify the uncertainty in data.  These uncertain data are then used to build computational models for cancer risk predictions \cite{Crespo21}.} and environmental applications \cite{Khora19,Eck17,Carstens18,Crespo21,Wang24,Dong24}, and when the measurements are noisy \cite{Haftbaradaran08}. Another situation when overfitting must be prevented arises when the computational cost of evaluating the objective function or the requirements is high, thereby limiting the value of $n$. Insufficient and inaccurately measured data threaten the relevance of a data-driven design because such a design might be infeasible in practice or it might attain an objective value far greater than the one predicted. These anomalies will be mitigated by considering {\em perturbed} scenarios. In particular, we wish to prevent $\theta^\star$ from overfitting ${\mathcal D}$ by ensuring that the success domain $\mathcal{S}(\theta^\star)$ not only contains the scenarios therein but also the vicinity around them.

Second, we want to explore the trade-off between performance and robustness by solving (\ref{sta}), so the requirements are satisfied for only a fraction of the $n$ scenarios. By making this fraction as large as possible we obtain the most robust design and a comparatively large $J(\theta^\star)$. However, this value can be lowered by considering a smaller fraction.  The scenarios excluded from this fraction, called outliers\footnote{Note that this definition is rooted on the feasibility of a scenario instead of its relationship to other scenarios (an outlier is commonly defined as a point that does not follow the trends followed by most of the data).}, will be chosen optimally during the search for $\theta^\star$ by using risk-averse or risk-agnostic formulations. These formulations account for or ignore the loss resulting from eliminating outliers respectively. The elimination of outliers mitigates the detrimental effects that extreme observations often have on data-driven designs \cite{Archimbaud18,Liang20,Guan19}. 


The main contributions of this paper are as follows:
\begin{enumerate}
\item	 As opposed to RO, RBDO, and DRO, the scenario-based strategies proposed herein make direct use of the available data, thereby eliminating the need for modeling the aleatory uncertainty, 
\item	 As in DRO, these strategies take into account the effects of only having a limited amount of data, and of that data possibly being corrupted by measurement error and uncertainty. In contrast to DRO, however, the proposed formulation accounts for these effects by assigning a perturbational set to each datum,
\item In contrast to DRO, this approach is readily applicable to implicitly-defined response and requirement functions requiring numerical simulation for their evaluation.
\item	Designs that trade-off performance against robustness are computed by eliminating optimally chosen elements of both the training set and the perturbational set, 
\item	The computational effort of the design optimization process is reduced by using worst-case perturbations and by sequentially augmenting a small training set. The example below shows that a Monte Carlo analysis with 10000 samples of an aeroelastic wing trained with only 50 scenarios attains
a good performance without exhibiting any failure. 
\item	The generalization properties of the resulting designs, both nominal and perturbational, can be evaluated by using Monte Carlo analysis and Scenario theory. The latter analysis is a rigorous, non-asymptotic, distribution-free assessment that accounts for the infinitely many elements in the perturbational sets, i.e., it accounts for both the ambiguity in the distribution and the uncertainty/error in the measurements leading to the data.
\end{enumerate}

This paper is organized as follows. Section \ref{UM} introduces a few uncertainty model types for the perturbed scenarios. This is followed by Section \ref{moti}, where key concepts are motivated using an engineering example. Several risk-averse and risk-agnostic formulations to robust design are presented and exemplified in Sections \ref{momi} and \ref{MOM}. This is followed by Section \ref{fea}, where we evaluate the reliability and robustness of the resulting designs. Strategies aiming to lower the computational cost of computing the proposed data-driven designs are introduced in Section \ref{lowcost}. 
In Section \ref{wing} we use the proposed framework to design several aeroelastic wings.  Finally, a few concluding remarks close the paper. 

\section{Modeling Uncertainty in the Data} \label{UM}
The pursuit of a design that remains feasible when the data are perturbed, a property called {\em robustness} hereafter, entails modeling such perturbations. A model for such perturbations should be prescribed according to the accuracy and precision of the metrology system, and/or to the analyst's belief of where the true value of $\delta$ might be.  While overly small perturbations lead to designs that violate the requirements in practice, overly large perturbations lead to conservative designs having an unnecessarily high $J(\theta^\star)$.

One model class is given by compact sets. Denote $\delta^{(i)}_s=\{\delta: \|\delta-\delta^{(i)}\| \leq \mu^{(i)}\}$ as the model of a perturbed scenario, where $\|\cdot\|$ is a norm and $\mu^{(i)}\geq 0$ is given by a rule prescribed by the analyst. Different norms yield to different set geometries. For instance, the $\ell^2$-norm makes $\delta^{(i)}_s$ a sphere centered at $\delta^{(i)}$ of radius $\mu^{(i)}$. As such, $\delta^{(i)}_s$ is fully prescribed by a realization of the unknown distribution $\mathbb{P}_\delta$ and the above mentioned rule. This model class leads to the sequence 
\begin{equation}
{\mathcal D}_s=\left \{\delta^{(i)}_s\right \}_{i=1}^n.\label{dsa}
\end{equation}

Another model class is given by distributions. Denote $\delta^{(i)}_d$  as a random vector having $\delta^{(i)}_s$ as the support set.  Hence, the uncertainty model $\delta^{(i)}_d$ is prescribed by a realization of the unknown distribution $\mathbb{P}_\delta$, and a distribution $\mathbb{P}_{\delta^{(i)}}$ set by the analyst. For instance, $\delta^{(i)}_d$ might be a Normal distribution having  $\delta^{(i)}$ as its mean and a fixed covariance matrix. This model class leads to the sequence 
\begin{equation}
{\mathcal D}_d=\left \{\delta^{(i)}_d\right \}_{i=1}^n.\label{dsa2}
\end{equation}
 
Designs that are robust to uncertainty in $\delta$ can be pursued by various means. The strategies below make use of a multi-point representation of (\ref{dsa2}). Let $\delta^{(i)}_p=\{\delta^{(i,j)}\}_{j=1}^m$ be a collection of $m$ sample points drawn from $\delta^{(i)}_d$. 
This leads to the multi-point sequence 
\begin{equation}
{\mathcal D}_p(m)=\left \{\delta^{(i)}_p \right \}_{i=1}^n.\label{dsc}
\end{equation}
An element of (\ref{dsa}), (\ref{dsa2}) or (\ref{dsc}) will be called a {\em perturbed} scenario. 

The data sequence in (\ref{dsc}) lead to the following classification: (i) design points computed from a single scenario will be called {\em single-point}, (ii) design points computed from more than one scenario will be called {\em multi-point}, and (iii) multi-point designs for which $m>1$ will be called {\em multi-point-robust}. 

The forthcoming programs are reformulations of 
\begin{align}
	 \min_{\theta \in \Theta} & \quad J\left(\theta \right) \label{sta2}\\
	\text{ subject to:} & \quad r_k(\theta,\delta^{(i,j)})\leq 0, \text{ for } k=1,\ldots n_r  \text{ and most }\delta^{(i,j)}\in{\mathcal D}_p. \nonumber 
\end{align} 
The elements of ${\mathcal D}_p$ for which $r_k(\theta^\star,\delta^{(i,j)})> 0$ for an unacceptably large number of the $m$ points in $\delta_p^{(i)}$ will be called outliers. Hereafter, the set of outliers will be denoted as ${\mathcal O}$ with $|{\mathcal O}|=\sigma$, whereas its complement set, $I$, will denote the set of inliers. The selection of outliers, which is made optimally as $\theta^\star$ is searched for, might depend on the severity of the corresponding requirement violation. This feature is explained next. Risk is commonly defined as the probability of an adverse outcome times the corresponding ``loss'' or consequence. In the context of this paper, the adverse outcome is the violation of a requirement, whereas the loss measure is the corresponding positive value taken by $r(\theta,\delta)$. A design $\theta^\star$ that depends on a loss measure will be called \emph{risk-averse}. This dependency limits the extent by which the outliers are allowed to violate the requirements. In contrast, formulations in which $\theta^\star$ does not depend on a loss measure will be called {\em risk-agnostic}.  

As expected, multi-point designs might not exhibit the desired degree of robustness,  e.g., most sets in ${\mathcal D}_s$ might not be contained in $\mathcal{S}(\theta^\star)$ even though all the points in ${\mathcal D}_p$ are. The robustness of $\theta^\star$ can be increased by using greater values for $\ell^{(i)}$, $n$ and $m$ in (\ref{sta2}). Unfortunately, the semi-infinite programming approaches that eliminate this error are inapplicable to general requirement functions. In the case when they are applicable, e.g., the functions are polynomial and the perturbed sets are bounded and semi-algebraic, the high computational cost limits their usage to moderate $n_\delta$ values \cite{Lacerda17}.  As such, (\ref{sta2}) is a computationally viable heuristic to (\ref{sta}).

\section{Motivational Example}\label{moti}
Consider the design of an aeroelastic wing subject to a structural dynamic constraint. In this setting, the decision variable $\theta$ parameterizes the wing's geometry, whereas the parameter $\delta$ is comprised of structural damping coefficients, which are intrinsically uncertain, and the Mach number, which varies with the wing's operating condition. The objective function $J(\theta)$ is the mass of the wing, whereas $r(\theta,\delta)>0$ implies flutter instability. Therefore, (\ref{sta2}) seeks the lightest wing that is stable for (most of) the $n$ scenarios in 
${\mathcal D}$. 

Note that the underlying physics might preclude $\theta^\star$ from satisfying the system requirement for some scenarios. This outcome could be the result of a limiting wing topology, a poor materials selection, or of an extreme/unrealistic scenario. A multi-point wing design that satisfies the requirement for as many scenarios as possible can be obtained by using the formulations in\footnote{To find the most robust risk-averse design use $J=0$. To find the most robust risk-agnostic design use $J=\|\alpha\|$ after making $\alpha$ an additional decision variable. These settings lead to designs that minimize the number of outliers.} Section \ref{momi}.  These formulations not only render $\theta_1^\star$ but also the corresponding list of $\sigma_{1}$ outliers. If $\sigma_{1}$ is unacceptably large, new wing topologies and different materials should be considered. Assume that $\sigma_{1}$ is acceptably low. At this point the analyst might wonder how much lighter the wing could become if we let the wing flutter for a few more scenarios. The same formulations can be used to derive a design having $\sigma_{2}>\sigma_{1}$ outliers. If $\theta_2^\star$ denotes such a design, it might turn out that $J(\theta_2^\star)$ is considerably smaller than $J(\theta_1^\star)$. At this point it is well worth pondering if the mass reduction $J(\theta_1^\star)-J(\theta_2^\star)$ justifies the loss in robustness caused by increasing the number of outliers from $\sigma_{1}$ to $\sigma_{2}$. 

Imagine now that if the $n$ scenarios in ${\mathcal D}$ are allowed to vary uniformly less than 5\% from their nominal value, the wing design $\theta_2^\star$ will be unstable for $\sigma'_{2}>\sigma_{2}$ of them. High sensitivity of $\sigma'_{2}-\sigma_{2}$ with respect to small variations in $\delta$ indicates that $\theta_2^\star$ is overfitting of the training set ${\mathcal D}$.  The poor generalization properties associated with this outcome might render $\theta_2^\star$ unacceptable.  To mend for this defficiency, a multi-point-robust wing is computed from the formulations in Sections \ref{FSR} or \ref{SCC1} by using the training set ${\mathcal D}_p(m)$. Say, the resulting wing $\theta_3^\star$ satisfies the requirements for as many perturbed scenarios as possible with an acceptably large probability, $\gamma$, set by the analyst. 
As before, the physics might prevent $\sigma_{3}$ from being zero. 
Designs having a lower mass and less robustness to uncertainty can be sought by lowering $\gamma$ and/or by accepting a greater number of outliers. 

Either physics-based limitations or choice might make a scenario an outlier. The process by which outliers are selected within the proposed formulations might depend on a risk (or loss) measure (Section \ref{FSR}) or not (Sections \ref{SCC1} or \ref{SCC2}). This measure quantifies the consequence of an adverse outcome, which in the context of this paper is the positive value of $r(\delta, \theta)$ attained by an outlier. 
For instance, a risk-averse design results from using $J(\theta)+{\mathbb E}[r(\delta^{(i,j)},\theta)\,|\,r(\delta^{(i,j)},\theta)> 0]$ as the objective function, where ${\mathbb E}[\cdot]$ is the sample mean. Note that the penalty term measures the strength of the instability. Conversely, a risk-agnostic approach ignores the possibly large values that this measure might take. For instance, a risk-averse design results from removing a few scenarios from ${\mathcal D}_p$ before solving (\ref{sta2}). In this setting, the strength of the instability of such scenarios is immaterial to the resulting $\theta^\star$. In contrast to this practice, the strategies below carry out an optimal selection of outliers, thereby leading to lower objective values.
Risk-averse and risk-agnostic approaches might lead to designs having a comparable probability of failure but very different objective values and loss measures. 

Furthermore, we might consider a wing design in which material densities are also uncertain, thereby making (\ref{sta2}) depend on moments of a response function. The developments in Section \ref{MOM} enable seeking the wing of minimal expected mass after the same outliers are eliminated from both the objective and the requirement functions. The final design should be chosen from several design alternatives having varying degrees of robustness and performance after comparing their objective function value, number of outliers, and loss measures.  A realistic wing design example is presented below. 
The formulations leading to various wing designs are detailed next. 
 
\section{Optimization Programs without Moments}\label{momi}
In this section we present risk-averse and risk-agnostic reformulations of (\ref{sta2}) in which neither the objective function nor the requirement functions depend on moments.  The resulting CCed programs use approximations to the {\em Cumulative Distribution Function} (CDF) of a random variable $z$ and its inverse based on the function evaluations in ${\mathcal Z}=\{z^{(i)}\}_{i=1}^n$. These approximations, denoted as $F_{\mathcal Z}$ and $F^{-1}_{\mathcal Z}$ respectively hereafter, are detailed in the Appendix.

\subsection{Risk-averse Scenario-based Formulation}\label{FSR}
Consider the optimization program
\begin{align}
	\label{chancon_sce}
	\min_{\theta\in \Theta,\; \xi \geq0} & \quad J\left(\theta\right)+\rho\sum_{i=1}^n \xi_{i}\\
	\text{subject to:} 
	                           & \quad F^{-1}_{{\mathcal N}(\theta,\,\delta_p^{(i)},\,k)}(\gamma_k)\leq \xi_{i},\;  i=1,\ldots n, \; k=1,\ldots n_r, \nonumber
\end{align}
where $\rho\in \mathbb{R}$ is a penalty parameter, $\xi \in \mathbb{R}^n$ is a slack variable, 
\begin{equation}
\label{za}
{\mathcal N}(\theta,\delta_p^{(i)},k)\triangleq \Big\{ r_k(\theta,\delta^{(i,\,j)})\Big \}_{j=1}^m,
\end{equation}
is the sequence of function evaluations of the $k$th requirement for the $m$ sample points of the $i$th scenario, and $0\ll \gamma_k \leq 1$ is the minimally acceptable probability of success\footnote{The CC in (\ref{chancon_sce}) is equivalent to $\mathbb{P}_\delta[r_k(\theta,\delta)\leq \xi_i]=F_{{\mathcal N}(\theta,\,\delta_p^{(i)},\,k)}(\xi_{i})\geq \gamma_k$. Note that the value taken by the CDF might plateau at zero or one as $\theta$ is varied, thereby making gradient-based searches stop at a local extremum. As such, CCs are better cast in terms of the inverse CDF. }. Denote as $\theta^\star$ and $\xi^\star$ the solution to (\ref{chancon_sce}). Hence, (\ref{chancon_sce}) seeks a design that minimizes the sum of $J$ and a penalty term while bounding the individual success probabilities for most elements of (\ref{dsc}) from below. Note that a constraint of (\ref{chancon_sce}) corresponding to a fixed $i$ and a fixed $k$ for $\xi_i=0$ is a computationally tractable heuristic for $\mathbb{P}_{\delta^{(i)}}[{\mathcal S}_k(\theta)]\geq \gamma_k$. The CC in (\ref{chancon_sce}) is a reformulation of the constraints in (\ref{sta}), in which ``most'' refers to a subset of $\Delta$ where the $k$th requirement function is satisfied by the perturbed scenarios with probability no less than $\gamma_k$.

Outliers are the perturbed scenarios for which at least one of the individual probabilities of success is below the acceptable threshold:
\begin{equation}
{\mathcal O}=\left \{\delta_p^{(i)}\in{\mathcal D}_p: \;\max_{k=1,\ldots n_r}F^{-1}_{{\mathcal N}(\theta,\,\delta_p^{(i)},\,k)}(\gamma_k)>0 \text{ for }i=1,\ldots n \right \}. \label{outset}
\end{equation}
The constraint in (\ref{outset}) corresponding to $\theta^\star$ is equivalent to $\xi_i^\star>0$. When $ \gamma_k < 1$ a few points of $\delta^{(i)}_p$ might fall into the failure domain regardless of such a scenario being an outlier or not. 

The feasible set of (\ref{chancon_sce}) is 
\begin{equation}
\mathcal{X}_{\text{robust}}=\bigcap_{i\in A,\; j\in B}\; \mathcal{X}\left (\delta^{(i,j)}\right ), \label{fea2}
\end{equation}
where $\mathcal{X}$ was defined in (\ref{feaset}). Relaxation entails turning an optimization program into one with either looser constraints or fewer constraints, so the feasible set is enlarged and the objective value $J(\theta^\star)$ is potentially lower.  Program (\ref{chancon_sce}) carries out two types of relaxation. The first type reduces the elements of $A$ in (\ref{fea2}) by lowering $\rho$. In particular, (\ref{chancon_sce}) enables the analyst to obtain various solutions as $\rho$ is varied from zero (no regret for having probabilities of success below $\gamma_k$) to infinite (infinite regret for having probabilities of success below $\gamma_k$). Large values of $\rho$ will drive all components of $\xi_{i}^\star$ to zero, thereby ensuring that as many success probabilities for the $k$th requirement as possible are no less than $\gamma_k$. Moderately large values of $\rho$ will make some $\xi_i^\star>0$, therefore yielding a design for which some probabilities are below $\gamma_k$. Hence, depending on the value of $\rho$, some perturbed scenarios are allowed to violate the requirements beyond acceptable limits for the purpose of lowering $J(\theta^\star)$.  However, this action has itself a cost, as expressed by the auxiliary variables $\xi_{i}$: if $\xi_{i} > 0$ the constraint $F^{-1}_{{\mathcal N}}(\gamma_k) \leq 0$ is relaxed to $F^{-1}_{{\mathcal N}}(\gamma_k)\leq \xi_{i}$ in exchange for a cost increase of $\rho\,\xi_{i}$. Note that the dependency of the number of outliers $\sigma$ on $\rho$ is implicit. The second relaxation type reduces the elements of $B$ in (\ref{fea2}) by setting the value of $\gamma_k$ for $k=1,\ldots n_r$ so requirement violations for a fixed fraction of the $m$ points in $\delta_p^{(i)}$ are allowed. 

Therefore, the feasible set of (\ref{chancon_sce}) is given by (\ref{fea2}) with $A\subseteq\{1,\ldots n\}$, and $B(i)\subset \{1,\ldots m\}$. The expansion of this set for an optimally chosen set of outliers and an optimally chosen fraction of the points in $\delta_p^{(i)}$ for $i=1,\ldots n$ might yield a considerably lower objective value without impacting most of the data, i.e., the probability of failure corresponding to most perturbed scenarios is not affected by these choices.  

Program (\ref{chancon_sce}) might not have a solution for which all $\xi_{i}^\star$ are zero, thereby making relaxation crucial. Relaxation enables the analyst to pursue the most robust design of a given architecture, as well as the identification of the scenarios for which the design specifications cannot be met. Notice that the number of decision variables in (\ref{chancon_sce}) grows with the number of scenarios, thereby limiting its applicability to moderately large datasets.  

The non-zero terms in the summation measure the extent by which the outliers violate the requirements, thereby making (\ref{chancon_sce}) risk-averse.  In this context, the adverse outcome is the violation of a requirement with an unacceptably large probability, whereas $\xi_i$ is the loss measure. Hence, the dependency of the penalty term on $\xi$ makes $\theta^\star$ a risk-averse design. 

The individual CCs in (\ref{chancon_sce}) can be written as $\mathbb{P}_\delta[ r_k(\theta,\delta)\leq \xi_i]\geq \gamma_k=1-\epsilon_k$. The feasible space of (\ref{chancon_sce}) is different from the feasible space of a program using the joint CC 
\begin{equation}
\mathbb{P}_\delta \left [\cup_{k=1}^{n_r} \{\delta: r_k(\theta,\delta)\leq \xi_i\}\right ]\geq \gamma=1-\epsilon, \label{jcc}
\end{equation}
 instead. The Bonferroni inequality, $\mathbb{P}_\delta[{\mathcal F}(\theta)]\leq \sum_{i=1}^{n_r} \mathbb{P}_\delta[ {\mathcal F}_i(\theta)]$, implies that Program (\ref{chancon_sce}) is a conservative approximation to the joint CCed program when $\epsilon > \sum_{i=1}^{n_r} \epsilon_i$. Conservative approximations have a smaller feasible set, so they might lead to a higher $J(\theta^\star)$.  Requirement functions exhibiting strong parameter dependencies and suboptimal choices for the $\epsilon_i$'s increase the conservatism in the approximation. Note however that (\ref{jcc}) is equivalent to the individual CC $\mathbb{P}_\delta[ r_{\text{max}}(\theta,\delta)\leq \xi_i]\geq \gamma$, where $r_{\text{max}}(\theta,\delta) \triangleq \max_{k=1,\ldots n_r} r_k(\theta,\delta)$ is the point-wise maximum of the requirement functions over $\theta$ and $\delta$. Therefore, a joint CC can be cast in terms of an individual CC without introducing any conservatism. This paper uses individual CCs because (i) they are simpler to implement, (ii) formulations based on the inverse CDF, easily computable in one dimension, lower the chance of the converged $\theta^\star$ being a local minima (Footnote 8), (iii) the minimally acceptable probability of success $\gamma_k$ might not be the same for all requirements, and (iv) the proposed framework can readily handle joint CCs by using non-conservative individual CCs.

\subsubsection{Worst-case formulation}\label{WC}
Program (\ref{chancon_sce}) with $\gamma_k=1$ for all $k=1,\ldots n_r$ is equivalent to 
\begin{align}
	\label{wc}
	\min_{\theta \in \Theta,\,  \xi \geq 0} & \quad J\left(\theta \right) + \rho \sum_{i=1}^n\xi_{i} \\
	\text{subject to:} & \quad r_k\left(\theta,\delta^{(i,j)}\right) \leq \xi_{i},\; i = 1, \ldots n,\; j=1,\ldots m,\; k=1,\ldots n_r.  \nonumber
\end{align}
This formulation \cite{Campi25} is called ``worst-case" because $\xi_i^\star$ depends on the element of $ \delta^{(i)}_p$ at which the requirement function takes on the greatest value. A few remarks on this particular case are discussed next. 

Note that the $m$ constraints corresponding to a fixed $i$ and a fixed $k$ for $\xi_i=0$ are a computationally tractable heuristic for $\delta_s^{(i)}\subset {\mathcal S}_k(\theta)$. As such, (\ref{wc}) is an approximation to the semi-infinite program with constraints $r_k(\theta, \delta)\leq \xi_i$ for all $\delta \in \delta_s^{(i)}$ with $i = 1, \ldots n$ and $k=1,\ldots n_r$. Furthermore, the elements of ${\mathcal D}_p$ having at least one parameter point in the failure domain will be the outliers. In contrast to (\ref{chancon_sce}), (\ref{wc}) is a convex program when the requirement functions are convex. 

The feasible set of (\ref{wc}) is given by (\ref{fea2}) with $A\subseteq\{1,\ldots n\}$ and $B=\{1,\ldots m\}$. Hence, relaxation is only attained by ignoring the requirements for the outliers.  As such, it only takes a single element of $\delta_p^{(i)}$ falling into the failure domain for (\ref{wc}) to penalize a requirement violation. 
The relaxation of $B$ in (\ref{chancon_sce}) might make the resulting designs differ considerably from those resulting from (\ref{wc}) for the same number of outliers. The possibly significant lower value $J(\theta^\star)$ of a CC design along with some risk tolerance might render it preferable unless safety is paramount \cite{Crespo24gnc}.

\subsection{Risk-Averse Requirement-based Formulation}\label{FRR}
The number of decision variables in the above formulations grows with the number of scenarios, rendering them impractical when the dataset is large. Furthermore, the loss measure in the penalty term is driven by the greatest values taken by the requirement functions, thereby possibly making a few requirements dominate the others. This problem can be avoided by scaling these functions, a difficult task due to their dependency on $\theta$.  The formulation below eliminates these drawbacks. 

Consider the CC program
\begin{align}
	\label{chancon_req}
	\min_{\theta\in \Theta,\;0\leq \zeta\leq 1} & \quad J\left(\theta\right)+\rho^\top \zeta\\
	\text{subject to:} & \quad F^{-1}_{{\mathcal Y}(\theta,\,{\mathcal D}_{p},\,k,\,\gamma_k)}(1-\zeta_k)\leq  0,\;  k=1,\ldots n_r, \nonumber
\end{align}
where $\rho\in \mathbb{R}^{n_r}$ is the penalty parameter, $\zeta\in \mathbb{R}^{n_r}$ is a slack variable,
\begin{equation}
\label{zz}
{\mathcal Y} (\theta,{\mathcal D}_{p},k,\gamma_k) \triangleq \left \{F^{-1}_{{\mathcal N}(\theta,\,\delta_p^{(i)},\,k)}(\gamma_k) \right\}_{i=1}^n,
\end{equation} 
is the sequence of the $\gamma_k$-quantiles for the $k$th requirement for all perturbed scenarios, and ${\mathcal N}$ is in (\ref{za}). Hence, (\ref{chancon_req}) seeks a design that minimizes the sum of $J$ and a penalty term while making the $1-\zeta_k$ quantile of ${\mathcal Y} (\theta,{\mathcal D}_{p},k,\gamma_k)$ for all $k=1,\ldots n_r$ fall into the success domain.  As before, this is a flexible scheme that enables the analyst to explore various solutions as $\|\rho\|$ is varied from zero to infinite. Large values of $\|\rho\|$ will make $\zeta_k^\star$ approach zero, thereby ensuring that as many individual probabilities of success as possible are acceptably large. Smaller values of $\|\rho\|$ will yield designs for which less individual probabilities of success are acceptable. The outliers of (\ref{chancon_req}) are also defined by (\ref{outset}).

Note that $\zeta_k^\star$ is the fraction of the scenarios for which the probability of success for the $k$th requirement is less than $\gamma_k$.  Further notice that making all the components of $\rho$ equal gives the same importance to all the requirements.  The relaxation of the CC by means of the tunable $\zeta$ prevents (\ref{chancon_req}) from becoming infeasible. Designs resulting from (\ref{chancon_sce}) are generally different from those based on (\ref{chancon_req})  because they use different loss measures. In particular, $\xi_i^\star$ takes on a value in the positive range of the requirement functions whereas $\zeta_k^\star$ takes on a value in $[0,1]$. More importantly, the number of decision variables and the number of constraints in (\ref{chancon_req}) depend on the number of requirements $n_r$ instead of the number of scenarios $n$ in (\ref{chancon_sce}), thereby having a lower the computational cost, e.g., the number of finite differences per design evaluation required by (\ref{chancon_req}) is $(n_r+1) \times (n_\theta+n_r)$ whereas that for (\ref{chancon_sce}) is $(n\,n_r+1) \times (n_\theta+n)$.
 

A key benefit of the free relaxation carried out by the above formulations is that they permit seeking designs that can only satisfy the requirements for a subset of the scenarios without having any prior information about this subset.  However, their risk-averse nature often renders a greater objective value. Conversely, the formulations that follow do not use a loss measure to carry out the relaxation, thereby becoming practically risk-agnostic. These formulations enable the analyst to prescribe upfront the desired number of outliers, i.e., the number of elements of $A$ in (\ref{fea2}), as well as the fraction of the $m$ perturbations allowed to violate the requirements,  i.e., the number of elements of $B$ in (\ref{fea2}). This property gives the name ``fixed" used hereafter. However, the particular set of outliers, i.e., the elements of $A$ and $B$ in (\ref{fea2}), is chosen by the optimization program. Note that setting an overly small number of outliers will make the optimization program infeasible.  Several strategies can be used to prevent this outcome\footnote{Infeasibility is avoided by first solving a risk-averse formulation with $\rho\gg 1$ and choosing $\alpha=\sigma/n$. Alternatively, the analyst can make $\alpha$ an additional decision variable in (\ref{chancon_sce2}) or (\ref{chancon_req2}), and use $J(\theta)=\|\alpha\|_2$ to first obtain $\alpha^\star$. Any future choice of $\alpha$ satisfying $\alpha\leq \alpha^\star$ will make such a program feasible.}. 

More importantly, fixed relaxations often yield lower objective values than their free counterparts. This is the result of the scenarios chosen as outliers falling in the upper quantiles of the distribution of the requirement functions, thereby yielding the greatest expansion of the feasible space.  These formulations also have the advantage of having a number of decision variables that does not increase with the number of scenarios or the number of requirements. This feature is critical to the non-convex CC programs commonly found in many engineering applications.  

\subsection{Risk-agnostic Scenario-based Formulation} \label{SCC1}
The risk-agnostic extension of (\ref{chancon_sce}) is 
\begin{align}
	\label{chancon_sce2}
	\min_{\theta\in \Theta} & \quad J\left(\theta\right)\\
	\text{subject to:} & \quad F^{-1}_{{\mathcal Z}(\theta,\,{\mathcal D}_{p})}(1-\alpha)\leq 0,  \nonumber 
	\end{align}
where 
\begin{equation}
{\mathcal Z}(\theta,{\mathcal D}_{p})\triangleq \left \{    \max_{k=1,\ldots n_r}F^{-1}_{{\mathcal N}(\theta,\,\delta_p^{(i)},\,k)}(\gamma_k)     \right\}_{i=1}^n, \label{zbar}
\end{equation}
is the sequence of the greatest quantiles from all requirements for all scenarios, 
and $0\leq \alpha \ll 1$ is the fraction of the $n$ scenarios to be considered as outliers.  Hence, $\theta^\star$ minimizes $J$ while ensuring that $\ceil{100(1-\alpha)}\%$ of the perturbed scenarios satisfy the requirements with an admissibly large probability of success.  

\subsection{Risk-agnostic Requirement-based Formulation} \label{SCC2}
The risk-agnostic extension extension of (\ref{chancon_req}) is 
\begin{align}
	\label{chancon_req2}
	\min_{\theta\in \Theta} & \quad J\left(\theta\right) \\
	\text{subject to:} & \quad F^{-1}_{{\mathcal Y}(\theta,\,{\mathcal D}_{p},\,k,\,\gamma_k)}(1-\alpha_k)\geq  0,\; k=1,\ldots n_r, \nonumber
\end{align}
where ${\mathcal Y}$ is in (\ref{zz}), and $0\leq \alpha_k\ll 1$ is the fraction of the $n$ scenarios violating the CC for the $k$th requirement. Hence, $\theta^\star$ minimizes $J$ while ensuring that  $\ceil{100(1-\alpha_k)}\%$ of the $n$ probabilities of success for all $k=1,\ldots n_r$ requirements are acceptably large. Note that the outliers violating a requirement might be different from the outliers violating another requirement, and $\ceil{\max_k \alpha_k}\leq \sigma/n\leq\ceil{\sum{\alpha_k}}$. 

Designs based on (\ref{chancon_sce2}) and (\ref{chancon_req2}) are practically risk-agnostic because $\theta^\star$ will not depend on more than a single outlier, i.e., the scenario at which ${\mathcal Z}$ takes the positive value closest to zero. This is a consequence of the number of scenarios prescribing the value of the inverse CDF in (\ref{Finvlin}), which is equal to one when $\alpha$ is a multiple of $1/(n-1)$ or two otherwise.

\vspace{0.2cm}
\noindent
{\bf Example 1 (Data Enclosure Design)}: Next we use a toy problem in $n_\delta=2$ dimensions to illustrate key features of the above formulations. In particular, we seek a set of minimal volume that encloses a set of points in ${\mathcal D}$.  These points are random draws of an unknown data generating mechanism. The requirement functions defining this set are
\begin{align}
r_1(\theta,\delta)&=\|\delta-c_1\|_2-u_1 \leq 0, \label{con1}\\
r_2(\theta,\delta)&=u_2-\|\delta-c_2\|_2\leq 0, \label{con2}
\end{align}
where $c_1\in{\mathbb R}^{n_\delta}$,  $c_2\in{\mathbb R}^{n_\delta}$,  $u_1\in{\mathbb R}$ and $u_2\in{\mathbb R}$ are the decision variables, thus $\theta=[c_1,u_1,c_2,u_2]$. The decision space is given by $u_1>0$, $u_2>0$ and $\|c_2-c_1\|_2\leq u_1$. Hence, the success domain ${\mathcal S}(\theta)$ is comprised of the points inside the circle $C_1=\{\delta: r_1(\theta,\delta)\leq 0 \}$ and outside the circle $C_2=\{\delta: r_2(\theta,\delta)\leq 0 \}$ such that the center of $C_2$ is in $C_1$.  Our goal is to find the design point $\theta$ leading to a success domain of minimal volume enclosing (most of) the scenarios in ${\mathcal D}$ and their vicinity.  

This volume is computed using
\begin{equation}
J(\theta)=\text{Vol}(\Delta)\,F_{Z(\theta,\,{\mathcal U}')}(0),
\end{equation}
where $Z(\theta,\,{\mathcal U}')=\{ \max_{k=1,2} r_k(\theta, u^{(i)})\}_{i=1}^{n_u}$, and the elements of ${\mathcal U}'=\{u^{(i)}\}_{i=1}^{n_u}$ are samples uniformly distributed in $\Delta$. To start, we assume that ${\mathcal D}$ is comprised of $n=15$ nominal scenarios. 

\begin{table}[h!]
\scriptsize
\setlength{\tabcolsep}{2pt}
\caption{Figures of merit for several designs. The area of the success domain is evaluated using $n_u=\num{2e4}$ sample points.}
\centering
 \label{tab}
 \begin{tabular}{|c| c c c c c c c c c c|}
 \hline 
    &  $\theta_1^\star$  & $\theta_2^\star$ &  $\theta_3^\star$  & $\theta_4^\star$ & $\theta_5^\star$  &$\theta_6^\star$  &$\theta_7^\star$  &$\theta_8^\star$  &$\theta_9^\star$  &$\theta_{10}^\star$      \\ [0.5ex]
\hline
\text{Formulation}                                        &  (\ref{wc}) &   (\ref{wc})  & (\ref{chancon_sce2})  &(\ref{chancon_sce2})  &(\ref{wc})   &(\ref{wc})    &(\ref{chancon_sce})  &(\ref{chancon_sce})    &(\ref{chancon_req2})          &(\ref{chancon_req2})     \\[0.5ex] 
$m$                                                            &  1             &   1        &1 &1         &81               &81              &81                    &81                     &81                              &81                        \\[0.5ex] 
$\sigma$                                                         &  0             &   1         &1 &2        &0                &1                &0                       &1                       &1                                &2                          \\[0.5ex] 
$J(\theta^\star)$                                          &  11.26     &   10.53   & 8.30 & 7.79    &15.13          &14.97        &14.63               &13.41                 &9.64                          &6.57                 \\[0.5ex] 
\hline
\end{tabular}
\end{table}

\begin{figure*}[h!]
    \centering
    \begin{subfigure}[b]{0.325\textwidth}
        \centering
        \centerline{\includegraphics[trim={1cm 6cm 3.5cm 6cm},clip,width=2.4\textwidth]{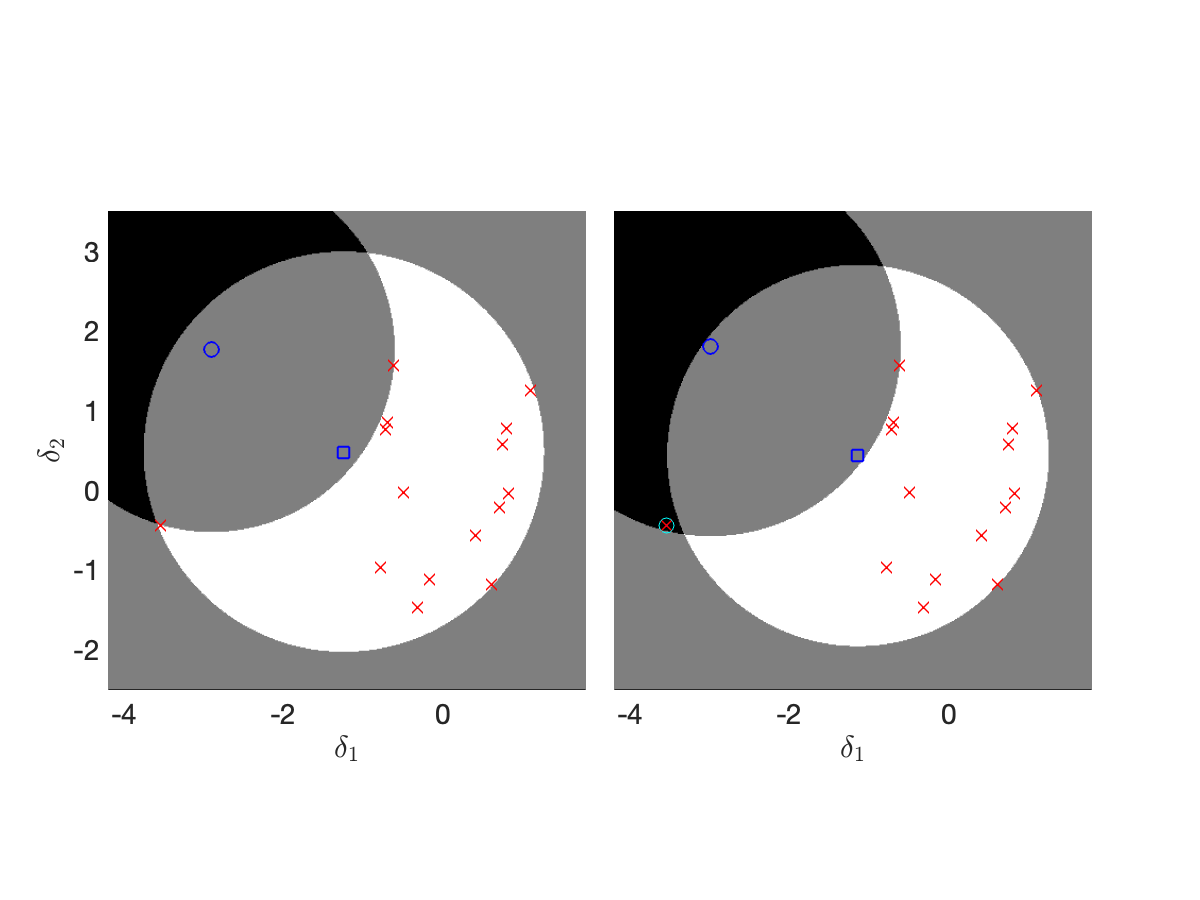}}
    \end{subfigure}
    \vfill
    \begin{subfigure}[b]{0.325\textwidth}
        \centering
        \centerline{\includegraphics[trim={0.6cm 7cm 3.5cm 9.5cm},clip,width=2.43\textwidth]{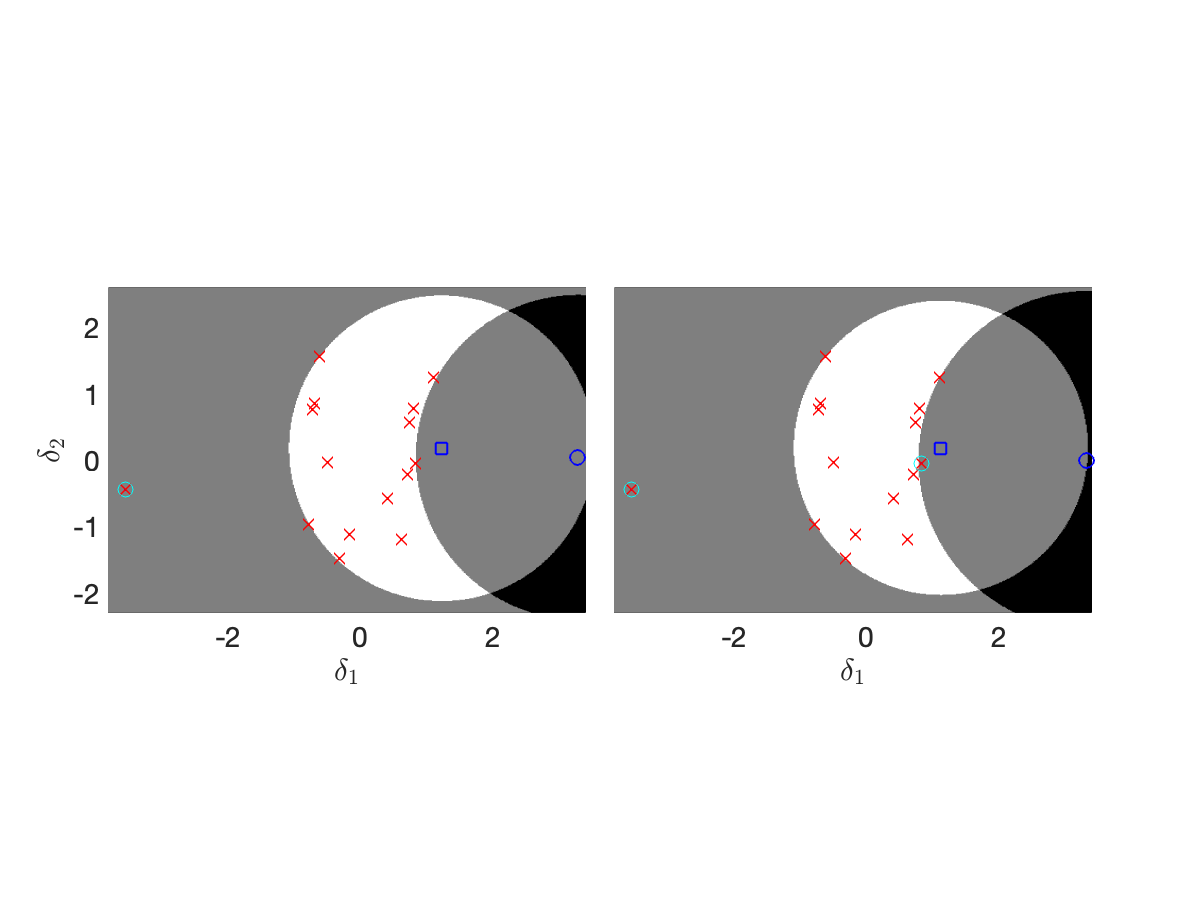}}
    \end{subfigure}
    \hfill
    \caption{\footnotesize{Success (white) and failure (non-white) domains for $\theta_1^\star$ (top-left), $\theta_2^\star$ (top-right), $\theta_3^\star$ (bottom-left) and $\theta_4^\star$ (bottom right).  The failure domain is shown in tones of gray: the darker the color the greater the number of constraint violations.  The nominal scenarios are marked with a ``{\color{red} $\times$}'' whereas the nominal scenarios falling into the failure domain are also marked with ``{\color{cyan} $\circ$}'' (if any). The center $c_1$ is marked with a ``{\color{blue} $\tinysquare$}'' whereas the center $c_2$ is marked with a ``{\color{blue} $\circ$}''.}}
    \label{fig1234}
\end{figure*}

Two multi-point designs, denoted as $\theta_1^\star$ and $\theta_2^\star$, derived from the worst-case formulation (\ref{wc}) for $m=1$ are presented first. The top subplots in Figure \ref{fig1234} show their corresponding success and failure domains along with the nominal scenarios. Whereas ${\mathcal S}(\theta_1^\star)$ encloses all the scenarios, ${\mathcal S}(\theta_2^\star)$ excludes an outlier. Table 1 lists relevant metrics for these and other designs. Because the formulation is risk-averse, the elimination of the outlier rendered a small reduction in $J(\theta^\star)$. This is the result of the penalty term in the objective function keeping the positive element of $\xi^\star$ near zero. Note that $\theta_1^\star$ is not robust to perturbations in the data. For instance, most perturbations of the nominal scenario $\hat{\delta}=[-3.7,-0.4]$ fall onto the failure domain. 

Two multi-point designs, denoted as $\theta_3^\star$ and $\theta_4^\star$, derived from the risk-agnostic formulation (\ref{chancon_sce2}) for $m=1$ are presented next. The bottom subplots of Figure \ref{fig1234} show their corresponding success and failure domains. Whereas ${\mathcal S}(\theta_3^\star)$ excludes a single outlier, ${\mathcal S}(\theta_4^\star)$ excludes two outliers.  Because this formulation is risk-agnostic, the elimination of the first outlier rendered a greater reduction in the objective value than $\theta_2^\star$. The elimination of the second outlier led to a slight improvement. However, all four designs are not robust to perturbations in the training data.  

\begin{figure*}[h!]
    \centering
    \begin{subfigure}[b]{0.325\textwidth}
        \centering
        \centerline{\includegraphics[trim={1cm 7cm 3.5cm 6cm},clip,width=2.4\textwidth]{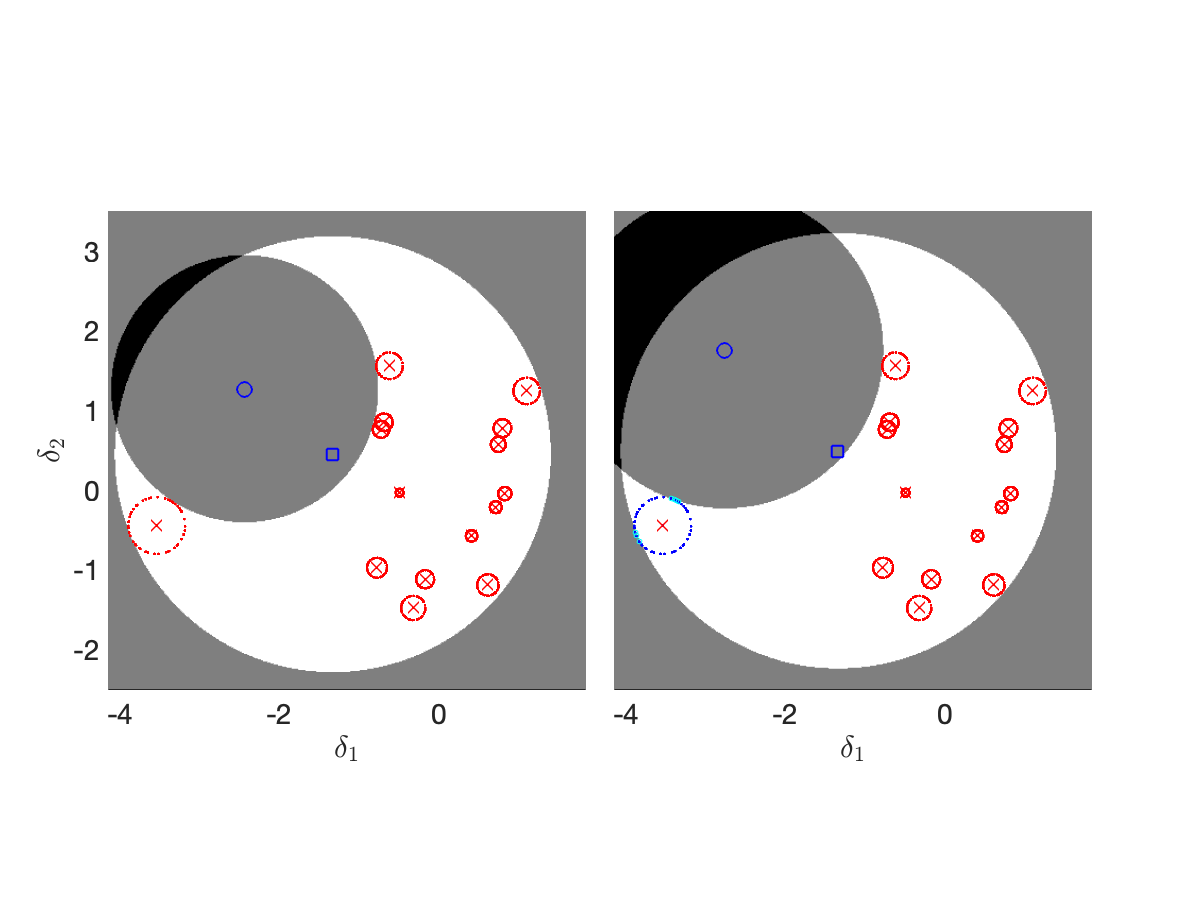}}
    \end{subfigure}
    \vfill
    \begin{subfigure}[b]{0.325\textwidth}
        \centering
        \centerline{\includegraphics[trim={0.98cm 6cm 3.5cm 7cm},clip,width=2.405\textwidth]{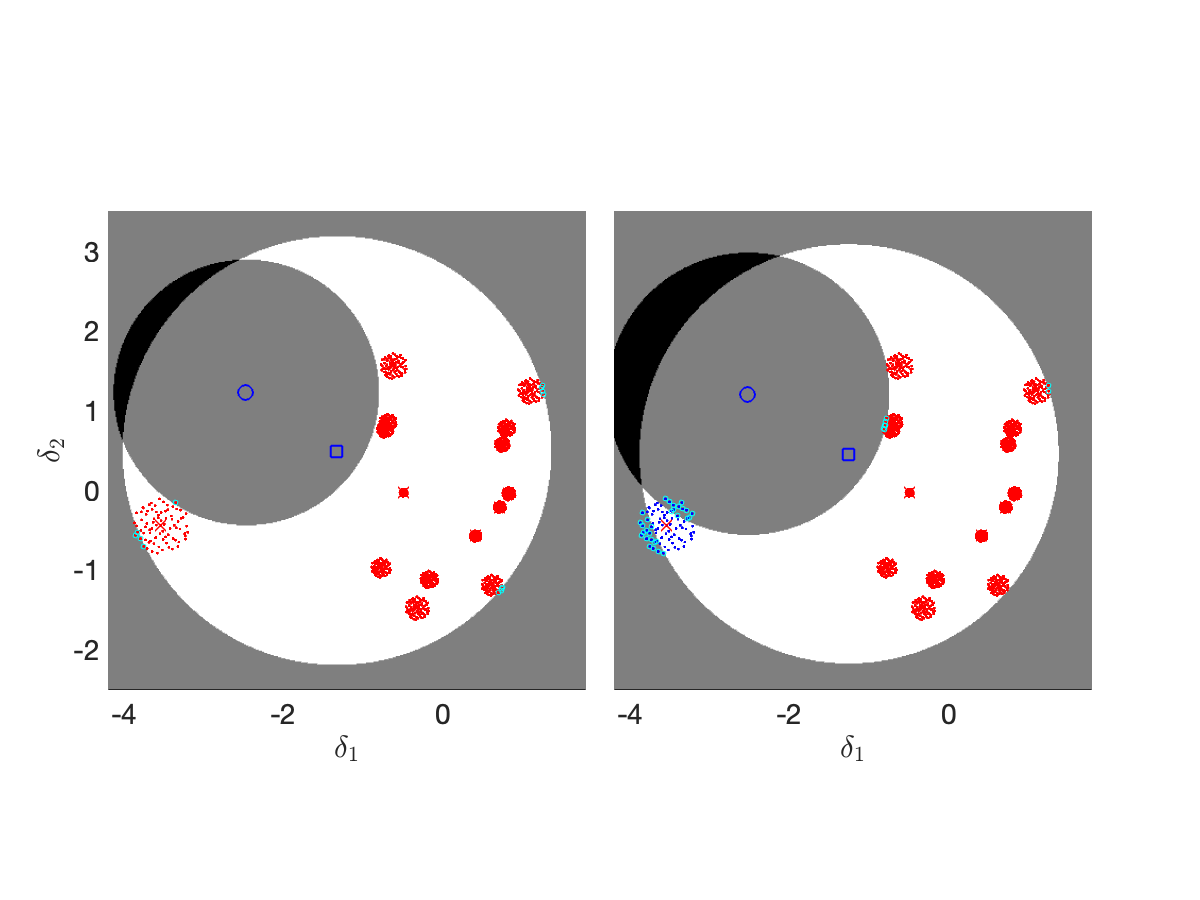}}
    \end{subfigure}
    \vfill
    \begin{subfigure}[b]{0.325\textwidth}
        \centering
        \centerline{\includegraphics[trim={0.6cm 7cm 3.5cm 9.75cm},clip,width=2.43\textwidth]{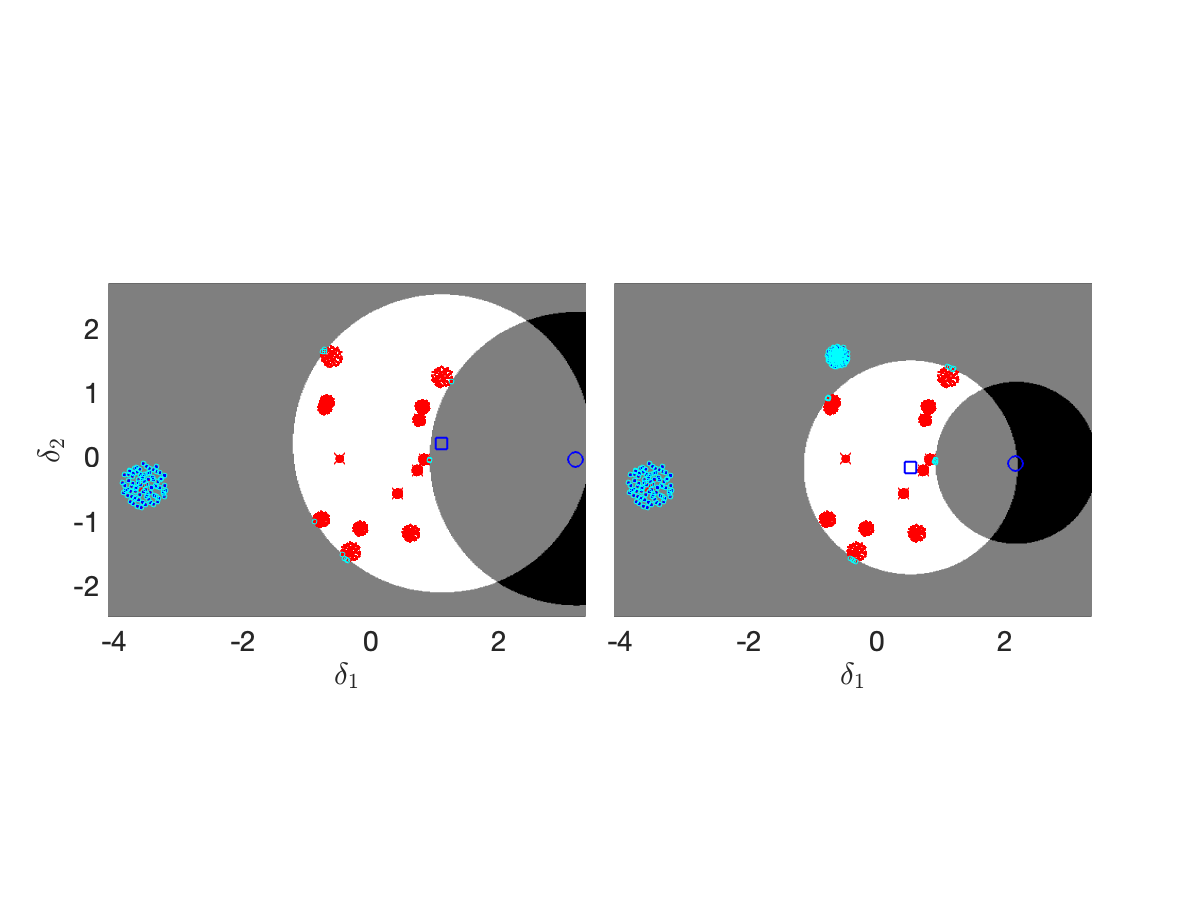}}
    \end{subfigure}
    \hfill
    \caption{\footnotesize{Success and failure domains for $\theta_5^\star$ (top-left), $\theta_6^\star$ (top-right), $\theta_7^\star$ (middle-left), $\theta_8^\star$ (middle-right), $\theta_9^\star$ (bottom-left) and $\theta_{10}^\star$ (bottom-right). The elements of $\delta_p^{(i)}$ corresponding to outliers are marked with a ``{\color{blue} $\boldsymbol{\cdot}$}''  (if any), inliers are marked with a ``{\color{red} $\boldsymbol{\cdot}$}'', and sample points falling into the failure domain are marked with a ``{\color{cyan} $\circ$}'' (if any).}}\label{fig3to8}
\end{figure*}

A few multi-point-robust designs are presented below. These designs are based on sequences having $m=81$ points distributed over circular support sets whose radii are proportional to the distance from the nominal scenario to the origin. Two sample sets ${\mathcal D}_{p}$ will be considered. The first sample set, shown in the top subplots of Figure \ref{fig3to8}, have points falling on the surface of the circles. This set is suitable for the worst-case formulation since the success domain is connected. The second sample set, shown in the middle and bottom subplots of Figure \ref{fig3to8}, have points falling on the area of the circles. This sample set is suitable for the CC formulations. 

Two worst-case designs derived from (\ref{wc}), denoted as $\theta_5^\star$ and $\theta_6^\star$ hereafter, are presented first.  The top subplots of Figure \ref{fig3to8} shows that ${\mathcal S}(\theta_5^\star)$ encloses all the perturbed scenarios, whereas ${\mathcal S}(\theta_6^\star)$ excludes an outlier.  As expected, the size of ${\mathcal S}(\theta_5^\star)$ is greater than that of ${\mathcal S}(\theta_1^\star)$, and the vicinity of $\hat{\delta}$ now falls into ${\mathcal S}(\theta_5^\star)$. The relaxation caused by eliminating an outlier made a few points of such a scenario fall into ${\mathcal F}(\theta_6^\star)$. However, the objective value only improved slightly. Designs $\theta_5^\star$ and $\theta_6^\star$ exhibit better robustness than their data-overfitting counterparts, but the loss measure in the penalty term prevents driving the scenario at $\hat{\delta}$ deeper into the failure domain, an outcome that would significantly lower $J(\theta^\star)$.

Two designs based on the risk-averse formulation (\ref{chancon_sce}) with $\gamma_1=\gamma_2=0.95$, denoted as $\theta_7^\star$ and $\theta_8^\star$ hereafter, are presented next. The parameter space corresponding to these designs are shown in the middle subplots of Figure \ref{fig3to8}. Whereas $\theta_7^\star$ drives all perturbed scenarios into the success domain with a probability not less than 0.95, $\theta_8^\star$ fails to do so for a single outlier. Note that a few points of three scenarios fall into ${\mathcal F}(\theta_7^\star)$ even though $\sigma=0$. $J(\theta^\star_8)<J(\theta^\star_6)$ thanks to the expansion of the feasible set attained by the CCs, i.e., $A$ in (\ref{fea2}) was kept the same whereas $B$ was reduced for a single scenario. By changing both $\rho$ and $\gamma$ the analyst can tune the robustness of a design.


Finally, the risk-agnostic formulation in (\ref{chancon_req2}) was used to derive $\theta_9^\star$ and $\theta_{10}^\star$. The bottom subplots of Figure \ref{fig3to8} show the corresponding parameter spaces. Whereas ${\mathcal S}(\theta_9^\star)$ contains all but one perturbed scenarios with a probability no less than 0.95, ${\mathcal S}(\theta_{10}^\star)$ excludes two scenarios. As before, a few points of 5 scenarios fall into $\mathcal {F}(\theta_9^\star)$ even though $\sigma=1$. That is also the case for $\mathcal {F}(\theta_{10}^\star)$, which contains points belonging to six perturbed scenarios. More importantly, note that both risk-agnostic optimal designs fully eliminate the effects of the outlier centered at $\hat{\delta}$, whose $m$ points violate the first requirement by a large margin. This leads to reductions in the objective value of about $37.5\%$ and $57.5\%$ relative to $\theta_5^\star$. This illustrates how the elimination of the worst-performing scenarios enforced by the risk-agnostic formulations yields a considerably lower objective value than the risk-averse formulations. This feature, along with having a number of decision variables that does not scale with the number of scenarios, might make risk-agnostic approaches preferable. Recall, however, that such approaches often lead to non-convex optimization programs even when the requirements are convex.

Naturally, the training dataset prescribes the resulting design. This dependency is stronger when $n$ is small. The developments in Section \ref{fea} evaluate the extent by which the properties enforced by design to the training dataset generalize to other datasets.

\section{Optimization Programs with Moments} \label{MOM}
In this section we consider problems that seek to minimize a moment of a response function. This is accomplished by using $J=\lambda$ and adding the constraint ${\mathbb M}[h(\theta,\delta)]<\lambda$ to (\ref{sta2}), where ${\mathbb M}$ is an empirical moment of the response function $h(\theta,\delta)$ and $\lambda$ is an additional decision variable. A key objective of the forthcoming formulations is to optimally identify and eliminate the same set of outliers from all the constraints. Thus, outliers will either yield large requirement function values or strongly contribute to ${\mathbb M}[h(\theta,\delta)]$. To simplify the presentation we will focus on the empirical mean, but extensions to other moments can be readily made. Denote as ${\mathbb E}\,[h(\theta,\delta), {\mathcal D}_p, {\mathcal W}]$ the empirical mean of $h(\theta,\delta)$ based on the scenarios in ${\mathcal D}_p$ and the weights in ${\mathcal W}=\{w^{(i)}\}_{i=1}^{n}$ with $w^{(i)}>0$, i.e.,
\begin{equation}
{\mathbb E}\,[ h(\theta,\delta), {\mathcal D}_p, {\mathcal W}]=\frac{1}{\sum w^{(i)} m}\sum_{i=1}^n \;\sum_{j=1}^m h(\theta,\delta^{(i,j)})\, w^{(i)}. \label{mea}
\end{equation}
The weights in $ {\mathcal W}$ will be used to select outliers. The programs below change the weight values during optimization. Once the optimization converges, the weights will take values near zero or one. Therefore, the weighted mean in (\ref{mea}) becomes the unweighted empirical mean for an optimally chosen subset of the data.

\subsection{Risk-averse Formulation}
Consider the optimization program
\begin{align}
	\label{mom_chancon_sce}
	\min_{\theta\in \Theta,\; \lambda,\; \xi\geq0} & \quad \lambda+\rho\sum_{i=1}^n \xi_{i}\\
	\text{subject to:} 	    & \quad F^{-1}_{{\mathcal N}(\theta,\,\delta_p^{(i)},\,k)}(\gamma_k)\leq \xi_{i},\;  i=1,\ldots n, \; k=1,\ldots n_r, \nonumber \\
	                                     & \quad {\mathbb E}\,[ h(\theta,\delta), {\mathcal D}_p, {\mathcal W}(\xi)] \leq \lambda, \nonumber 
\end{align}
where ${\mathcal W}(\xi)=\{\exp(-\kappa\, \xi_i)\}_{i=1}^{n_a}$ for $\kappa\geq1$ and ${\mathcal N}$ is in (\ref{za}). Denote as $\theta^\star$, $\lambda^\star$ and $\xi^\star$ the solution to (\ref{mom_chancon_sce}). Note that ${\mathcal W}(\xi^\star)$ takes the value of one when the scenario is an inlier, and the value of zero when the scenario is an outlier.  By making the weights depend smoothly on the decision variable $\xi$, gradient-based algorithms are applicable. Hence, (\ref{mom_chancon_sce}) seeks a design $\theta^\star$ that minimizes the sum of the empirical mean of the response for the inliers and a penalty term subject to the same CCs considered earlier. Note that the slack variable $\xi$ enables eliminating the same set of outliers from the moment-based and the non-moment-based constraints. In spite of this benefit, however, the designs resulting from (\ref{mom_chancon_sce}) exhibit the subpar objective value and high cost typical of risk-averse formulations. 

\subsection{Risk-agnostic Formulation} 
Note that the decision variable $\xi$ of (\ref{mom_chancon_sce}), used to consistently eliminate outliers from all the constraints, is unavailable. A risk-agnostic formulation based on (\ref{chancon_sce2}) is presented next.  
Consider the moment sequence
\begin{equation}
\label{zab}
{\mathcal U}(\theta,{\mathcal D}_{p}) \triangleq \Big \{\, {\mathbb E}\,[ h(\theta,\delta), {\mathcal D}_p, \left \{K_{ik}\right \}_{k=1}^n]\, \Big \}_{i=1}^n,
\end{equation}
where $K$ is the Kronecker delta. Hence, the $i$th element of ${\mathcal U}$ is an empirical mean of the response for all $m$ sample points of the perturbed scenario $\delta_p^{(i)}$. Furthermore, define ${\mathcal V}(i)$ as the subsequence of ${\mathcal U}$ 
whose elements do not exceed the $i$th-order statistic of ${\mathcal V}$, and $g \in \mathbb{R}^{n_r+1}$ as
\begin{align}
g(i,\theta,\lambda,{\mathcal D}_{p},\gamma)& \triangleq\left[\; {\mathcal V}(i)- \lambda,\; F^{-1}_{{\mathcal N}(\theta,\,\delta_p^{(i)},\,1)}(\gamma_1),\;\ldots,\; F^{-1}_{{\mathcal N}(\theta,\,\delta_p^{(i)},\,n_r)}(\gamma_{n_r})\;\right],
\end{align}
where $ i\in\{1,\ldots n\}$. Hence, a scenario $i^\star$ for which $g(i^\star,\theta,\lambda,{\mathcal D}_{p},\gamma)\leq 0$ satisfies two properties. First, the mean of the response for the lowest $n_{i^\star}$ elements of ${\mathcal V}$ does not exceed $\lambda$. Second, $i^\star$ is an inlier. This sets the stage for the optimization program
\begin{align}
	\label{chancon_sce4}
	\min_{\theta\in \Theta,\, \lambda} & \quad \lambda\\
	\text{subject to:} & \quad F^{-1}_{{\mathcal E}(\theta,\,\lambda,\,{\mathcal D}_{p},\,\gamma)}(1-\alpha)\leq 0, \nonumber 
\end{align}
where $0\leq \alpha\ll 1$ is the fraction of the $n$ scenarios considered as outliers, and
\begin{equation}
\label{zab2}
{\mathcal E}(\theta,\lambda,{\mathcal D}_{p},\gamma) \triangleq \left \{\max_{j=1,\ldots n_r+1}\; g_j(i,\theta,\lambda,{\mathcal D}_{p},\gamma) \right \}_{i=1}^n.
\end{equation}
Hence, (\ref{chancon_sce4}) yields a design $\theta^\star$ that minimizes the lowest sample moment of the response for $\ceil{100(1-\alpha)}\%$ of the scenarios while ensuring that they satisfy the requirements with an acceptably large probability of success.  The $\max$ operator in (\ref{zab2}) ensures that the same set of outliers are removed from all the constraints. As before, the risk-agnostic formulation (\ref{chancon_sce4}) often leads to a lower objective values than the risk-averse formulation (\ref{mom_chancon_sce}) for the same number of outliers.

\section{Reliability and Robustness Analyses} \label{fea}
The formulations above use data-driven heuristics for $\mathbb{P}_{\delta^{(i)}}[{\mathcal S}(\theta)]\geq \gamma$ based on the training dataset of perturbed scenarios ${\mathcal D}_p$. As expected, further analyses are required to determine if $\mathbb{P}_{\delta^{(i)}}[{\mathcal S}(\theta^\star)]\geq \gamma$ is satisfied for other data. The generalization properties of data-driven designs can be evaluated by using either Monte Carlo analysis or scenario theory. A Monte Carlo analysis estimates the probability of requirement violations for any $\theta$ regardless of the means by which such a design was obtained. However, the additional data required to compute such an estimate might not be available, e.g., when the nominal scenarios are obtained experimentally. In contrast, scenario theory yields a rigorous upper bound to this probability for scenario-based optimal designs without requiring additional data. 

In the developments that follow we will carry out two analyses. The first analysis, which neglects the uncertainty in the data, evaluates the probability of nominal scenarios falling into the failure domain, i.e., $\mathbb{P}_{\text{nom}}[{\mathcal F}(\theta)]\triangleq \mathbb{P}_{\delta}[\delta\in{\mathcal F}(\theta)]$. The second analysis, which accounts for uncertainty in the data, evaluates the probability of perturbed scenarios falling into the failure domain, i.e., $\mathbb{P}_{\text{per}}^\gamma[{\mathcal F}(\theta)]\triangleq\mathbb{P}_{\delta}[ \mathbb{P}_{\delta_d}[\delta_d(\delta)\in{\mathcal F}(\theta)]>1-\gamma]$.  This probability will be called the perturbational failure probability hereafter. 
These analyses will be referred to as \emph{reliability analysis} and \emph{robustness analysis} respectively.  

Failure probabilities can be readily evaluated using Monte Carlo. In the nominal case the process entails generating a large testing dataset of $n'$ nominal scenarios, and finding the fraction of them falling onto ${\mathcal F}(\theta)$. In the perturbational case, the distribution $\delta_d$ in (\ref{dsa2}) corresponding to each of the $n'$ scenarios is first simulated to obtain $m' \gg m$ sample points. The desired probability is the fraction of the $n'$ scenarios for which more than $\lceil m' \times (1-\gamma) \rceil$ points fall onto ${\mathcal F}(\theta)$. The sampling error caused by using a finite number of samples to make this determination should be quantified using confidence intervals.  

When the training dataset is \emph{Independent and Identically Distributed} (IID)\footnote{Note that this assumption is about the process by which data is drawn from a data-generating-mechanism, not about the dependencies among the components of $\delta$. This is a mild assumption often satisfied in practice.} both reliability and robustness analyses can be carried out using scenario theory. These analyses yield a non-asymptotic (it applies to designs trained with any number of scenarios), distribution-free (it applies to any data-generating mechanism) upper bound to the failure probability without modeling the underlying distribution. The developments in \cite{Garatti22} are applicable to nominal data, whereas those in \cite{Campi25} are applicable to perturbational data. The robustness analysis decouples the design process, which was carried out by using a finite number of perturbations, from the analysis process, which accounts for the infinitely many perturbations in each element of (\ref{dsa}). As such, this framework enables the analyst to rigorously evaluate the robustness of the optimal design against adversarial actions of any strength. For instance, we can evaluate the robustness of $\theta^\star$ against the rule leading to the $\delta_s$'s in (\ref{dsa}), or against any another rule.  


\vspace{0.2cm}
\noindent
{\bf Example 2 (Reliability Analyses of the Data-enclosing Sets)}: The designs presented in Example 1 are based on a small dataset thereby possibly exhibiting a large failure probability. Four designs based on a dataset having $n=500$ scenarios and $\gamma=0.95$ are presented next.  In particular, we synthesize (i) a multi-point-design with $\sigma=25$ outliers and $m=1$, (ii) a multi-point design without outliers and $m=1$, (iii) a multi-point-robust design with $\sigma=25$ outliers, and (iv) a multi-point-robust design without outliers. The resulting designs, denoted as $\theta_A^\star$, $\theta_B^\star$,  $\theta_C^\star$ and $\theta_D^\star$ respectively, are obtained from (\ref{chancon_sce2}).  In contrast to the perturbations ${\mathcal D}_p$ used previously, we will consider adversarial perturbations. In particular, the radius of the circular perturbational set for the $i$th scenario is 
\begin{equation}
r^{(i)}=r_{\text{max}}\,\exp\left(-\max_{k=1,\,2}\; r_k(\theta,\delta^{(i)})^2\right).
\end{equation}
Therefore, the closer the nominal scenario is to the boundary of the failure domain the stronger the perturbation.  Adversarial actions are an additional source of non-convexity. 

\begin{figure*}[h!]
    \centering
    \begin{subfigure}[b]{0.325\textwidth}
        \centering
        \centerline{\includegraphics[trim={1cm 7.3cm 3.5cm 6cm},clip,width=2.4\textwidth]{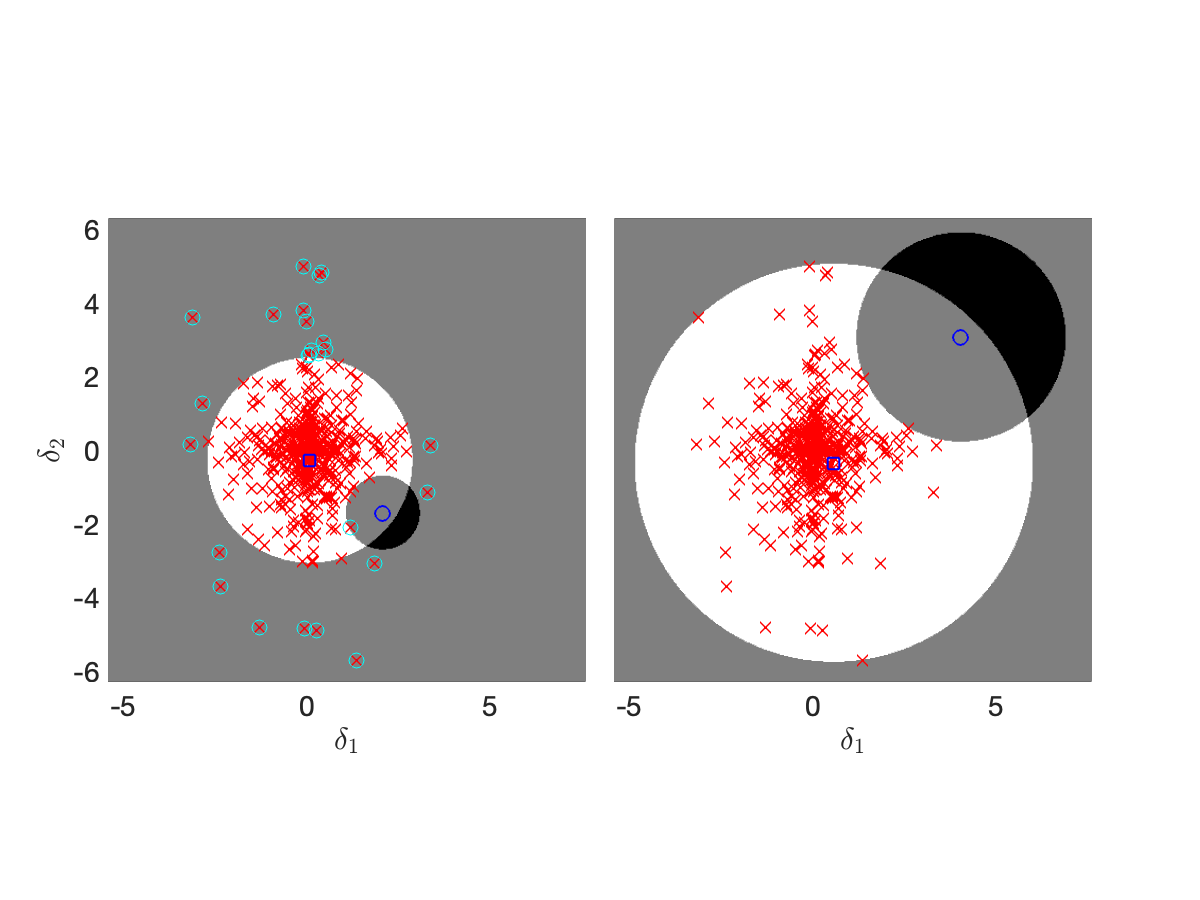}}
    \end{subfigure}
    \vfill
    \begin{subfigure}[b]{0.325\textwidth}
        \centering
        \centerline{\includegraphics[trim={0.6cm 5cm 3.5cm 7.3cm},clip,width=2.43\textwidth]{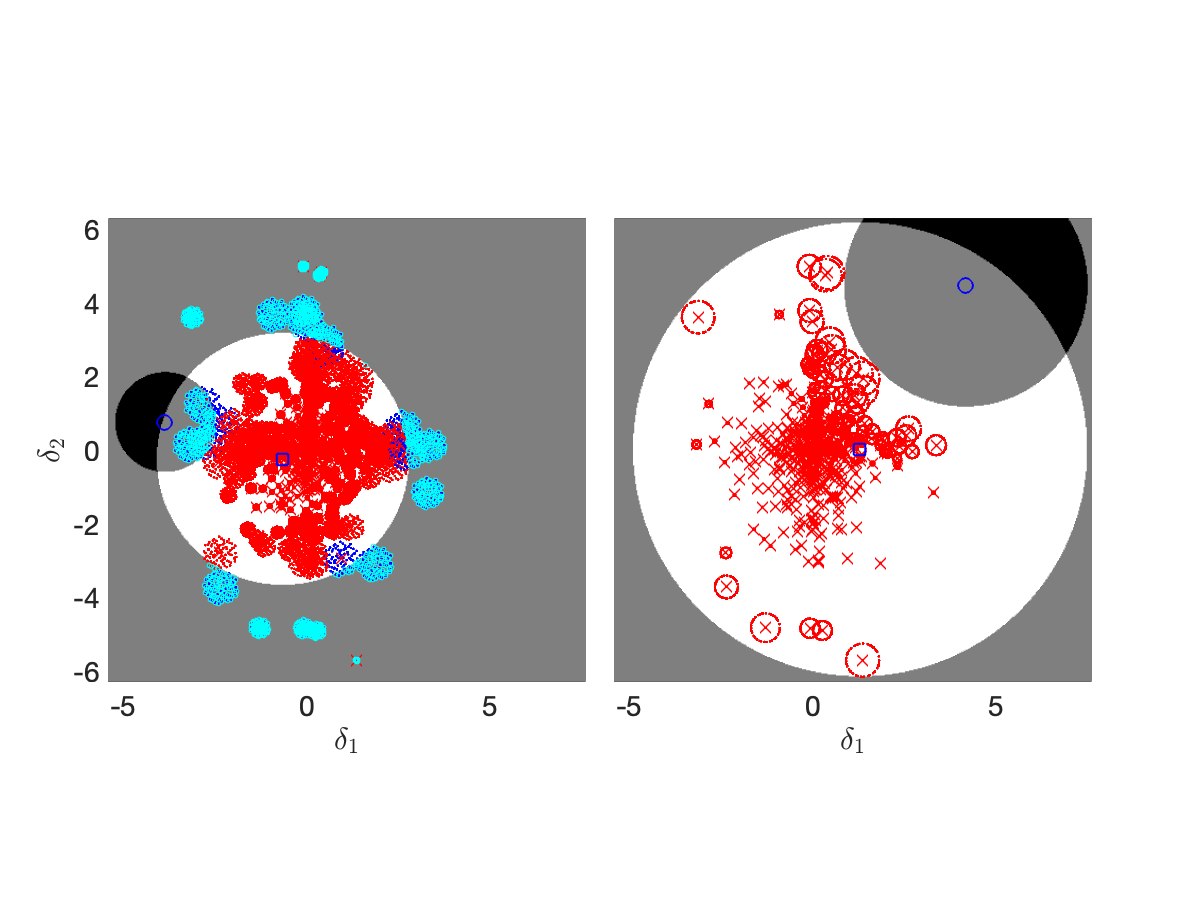}}
    \end{subfigure}
    \hfill
    \caption{\footnotesize{Success and failure domains for $\theta_A^\star$ (top-left), $\theta_B^\star$ (top-right),  $\theta_C^\star$ (bottom-left) and $\theta_D^\star$ (bottom-right).}}
    \label{figabcd}
\end{figure*}

Figure \ref{figabcd} shows the parameter spaces corresponding to the resulting designs. Such designs attain various degrees of performance and robustness depending upon the choices of $m$, $\sigma$, $\gamma$, and $r_{\text{max}}$.

\begin{table}[h!]
\scriptsize
\setlength{\tabcolsep}{2pt}
\caption{Performance and robustness metrics of several designs. $J(\theta^\star)$ is the objective value, $\mathbb{P}_{\text{nom}}[{\mathcal F}(\theta^\star)]$ is the nominal failure probability, $\mathbb{P}_{\text{per}}^\gamma[{\mathcal F}(\theta^\star)]$ is the perturbational failure probability, and CI is the $95\%$ confidence interval of these probabilities.
}
\centering
 \label{tab2}
 \begin{tabular}{|c| c c c | c c | c c|}
 \hline 
     & $m$ & $\sigma/n$  & $J(\theta^\star)$   &  $\mathbb{P}_{\text{nom}}[{\mathcal F}(\theta^\star)]$  & CI &  $\mathbb{P}_{\text{per}}^{\tiny 0.95}[{\mathcal F}(\theta^\star)]$ & CI  \\ [0.5ex]
\hline
$\theta_A^\star$                                          &     1  &  0.05     & 22.26   & \num{8.76e-2} & [\num{8.35},\,\num{9.17}]$\times 10^{-2}$ &  \num{1.46e-1}   &  [\num{1.40},\,\num{1.50}]$\times 10^{-1}$  \\[0.5ex] 
$\theta_{B}^\star$                                       &   1    & 0 & 77.68               & \num{1.60e-2}  & [\num{1.42},\,\num{1.78}]$\times 10^{-2}$ & \num{3.09e-2}           & [\num{2.87},\,\num{3.37}]$\times 10^{-2}$\\[0.5ex] 
$\theta_C^\star$                                          &   81 & 0.05      & 34.01     &\num{5.41e-2}  & [\num{5.09},\,\num{5.74}]$\times 10^{-2}$ & \num{8.34e-2}   & [\num{8.02},\,\num{8.83}]$\times 10^{-2}$\\[0.5ex] 
$\theta_{D}^\star$                                       &     81 & 0  & 91.36           &\num{8.95e-3}  &[\num{0.76},\,\num{1.03}]$\times 10^{-2}$ & \num{1.64e-2}           & [\num{1.47},\,\num{1.83}]$\times 10^{-2}$ \\[0.5ex] 
\hline
\end{tabular}
\end{table}

The reliability and robustness analyses of the four designs are presented in Table \ref{tab2}, where the objective value, the failure probabilities and their confidence intervals \cite{Hanson10} are provided. These estimates are based on a Monte Carlo campaign with $n'=20000$ nominal scenarios and $m'=200$ sample points. As expected, the perturbational probabilities are greater than the nominal probabilities. The analyst should choose a design among the alternatives according to the desired balance between performance, as measured by $J(\theta^\star)$, and robustness, as measured by the failure probabilities.  As expected, the performance-based ranking of the four designs is the opposite of the robustness-based ranking. Note that lowering the failure probability one order of magnitude increased the optimal objective value 4.13 times.  


\section{Lowering Computational Cost}\label{lowcost}
When the failure probability found in testing is unacceptably large, a new design should be computed using a greater $n$ and/or $m$.  A large training dataset is often required to accurately compute a small probability, thereby increasing the computational cost of solving the scenario program. 
Strategies that lower this cost without relaxing the underlying reliability/robustness specifications are outlined next.

\subsection{Adversarial Perturbations}\label{advper}
The pursuit of a multi-point-robust design requires prescribing perturbations of the nominal data. The greater the number of such perturbations, $m$, the higher the computational cost of solving for $\theta^\star$. This cost can be reduced by considering adversarial perturbations for only a subset of the scenarios.

In the context of this paper, an adversarial perturbation of the nominal scenario $\delta_n$ increases $r_k(\theta,\delta)$ the most near $\delta_n$. As such, these perturbations are locally worst-case in the $\delta$ continuum. A design that is robust to the worst-case perturbation of a set will be robust to all other perturbations therein. So, instead of using $m$ random perturbations to form $\delta^{(i)}_p$, we will use a single adversarial perturbation. This practice is beneficial when the cost of identifying this perturbation is lower than the cost of evaluating the requirement function $m$ times.  The adversarial perturbation of $\delta_n$ corresponding to the requirement $r_k$ for the baseline design $\hat{\theta}$ is  
\begin{equation}
\delta_a(\lambda,k,\delta_n,\hat{\theta})= \delta_n+\lambda \frac{\partial r_k(\theta,\delta)}{\partial \delta} \Big |_{\delta_n,\,\hat{\theta}}\,, \label{wcp}
\end{equation}
where $\lambda>0$. When multiple requirements are present we use $\delta_a(\hat{\lambda},\hat{k},\delta^{(i)},\hat{\theta})$, where $\hat{k}=\argmax_k  r_k(\hat{\theta}, \delta^{(i)})$ and $\|\delta_a(\hat{\lambda},\hat{k},\delta^{(i)},\hat{\theta})-\delta^{(i)}\|=\mu^{(i)}$ (the definition of $\delta_s^{(i)}$). Therefore, $\delta_a(\hat{\lambda},\hat{k},\delta^{(i)},\hat{\theta})$ is a point on the ball's surface. 

Further cost reductions are attained by only perturbing a small subset of the nominal scenarios. The scenarios to be perturbed must yield the greatest negative values in $\{r_{\hat{k}}(\hat{\theta},\delta^{(i)})\}_{i=1}^n$ (so the scenarios to be perturbed are the closest to the failure domain $\mathcal{F}_{\hat{k}}(\hat{\theta})$ as if $|r_{\hat{k}}(\hat{\theta},\delta^{(i)})|$ were a measure of distance).  Hence, the center of the ball $\delta_s^{(i)}$ is in the success domain, whereas $\delta_a(\hat{\lambda},\hat{k},\delta^{(i)},\hat{\theta})$ is either in $\mathcal{F}_{\hat{k}}(\hat{\theta})$ or close to $\mathcal{F}_{\hat{k}}(\hat{\theta})$.

In the above developments, the design $\hat{\theta}$ can be kept fixed, or it can be sequentially updated as the optimization algorithm converges to $\theta^\star$. The later option is computationally appealing when automatic differentiation is available.

\subsection{Sequential Design}\label{seque}
This approach entails generating a design sequence having increasing robustness levels based on comparatively small training sets. The basic idea is to sequentially augment the training dataset by adding a few testing points that marginally violate the requirements for the preceding design.  In this fashion, we leverage the computational cost of testing in order to systematically expand the success domain without using training scenarios not rendering robustness improvements.

For simplicity in the presentation we will focus on the case in which $m=1$ and the number of outliers is made as small as possible (Footnote 7). The same ideas can be extended to other settings. Denote as $\theta^\star_u({\mathcal D}_{p,u}(n_u))$ a robust design based on the training set ${\mathcal D}_{p,u}$ of $n_u$ scenarios, where the subindex $u\in \{1,2,\ldots\}$ is the iteration number. To start off, compute $\theta^\star_1({\mathcal D}_{p,1}(n_1))$ for a small value of $n_1$. Next, we perform a Monte Carlo analysis of $\theta^\star_1$ using the testing set ${\mathcal T}=\{\delta^{(i)}\}_{i=1}^{n'}$ with $n' \gg n$. If the probability of failure is unacceptably large, we identify a few sample points of ${\mathcal T}$ leading to positive but small $r_{\hat{k}}(\theta^\star_1,\delta^{(i)})$ values\footnote{When the testing scenarios are drawn from a synthetic distribution, such points must also attain a comparatively large likelihood.}. Hence, these points not only fall onto the failure domain $\mathcal{F}_{\hat{k}}(\theta^\star_1)$ but are also near its boundary. Denote as ${\mathcal A}=\{\delta^{(i)}\}_{i=1}^{n_a}$ the collection of such points. We will then compute  $\theta^\star_2({\mathcal D}_{p,2}(n_2))$, where ${\mathcal D}_{p,2}={\mathcal D}_{p,1}\cup {\mathcal A}$ and $n_2=n_1+n_a$. This process is repeated until the failure probability is either acceptable or the design sequence converges.

Note that the enlarging training sets are not comprised of IID observations. As such, the statistics used in training are not accurate representations of the statistics found in testing, e.g., the number of outliers in the last iteration might not accurately represent the probability of failure. Further notice that this feature invalidates the scenario-based bounds. 

\section{Robust Design of a Aeroelastic Wing} \label{wing}
Next we consider the design of a flexible aeroelastic wing subject to static and dynamic structural requirements. The objective function to be minimized is the sample mean of $h(\theta,\delta)$, which is a weighted combination of the structural wing mass, and the aerodynamic drag coefficient computed at a static aeroelastic trim condition. When a flexible wing is subjected to low dynamic pressures, structural dynamic perturbations will dampen out in time, thereby leading to an asymptotically stable response.  At higher dynamic pressures, however, coupling between the structure and the unsteady aerodynamics might result in an unstable wing.  This phenomenon is called aeroelastic flutter, and the first design requirement considered here, $r_1(\theta,\delta)<0$, ensures that the wing does not flutter below some prescribed dynamic pressure threshold.  The second design requirement, $r_2(\theta,\delta)<0$, ensures that the static stresses which develop along the wing at the trim condition do not exceed a limit. Hence, the success domain, ${\mathcal S}(\theta)$, is comprised of the $\delta$ points for which the wing does not flutter, and the peak stress is acceptably low. The design goal is to find the wing of minimal objective value that satisfies the requirements for most of the scenarios in ${\mathcal D}$ and their vicinity.

\begin{figure}[htb]
 \centering
 \includegraphics[trim={6cm 0cm 6cm 0cm},clip,width=2.5in]{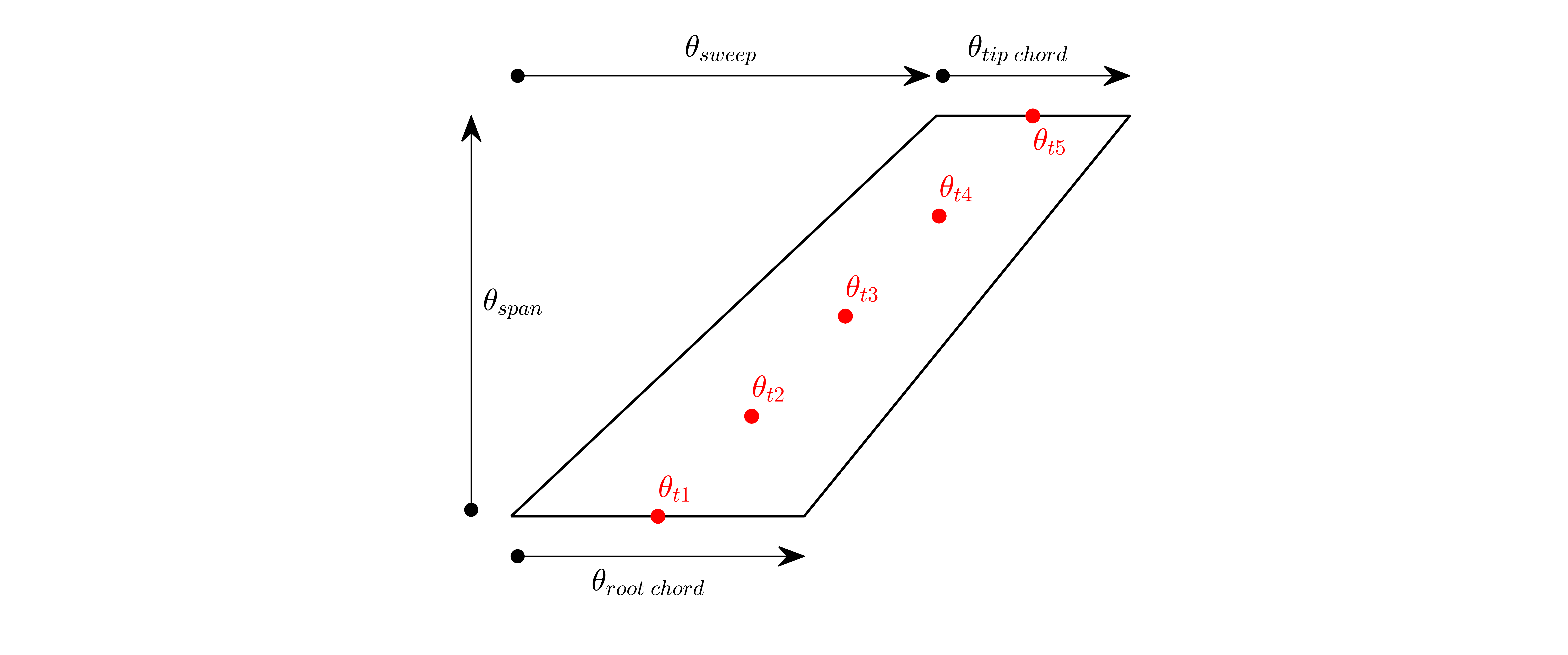}
 \caption{Shape and sizing design variables for the aeroelastic wing.}
 \label{wing_pic}
\end{figure}

The design variables $\theta\in  \mathbb{R}^{9}$, shown in Figure ~\ref{wing_pic}, are segregated into shape and sizing variables. The dimensions of the wing are set to emulate a subsonic wind tunnel test model \cite{Agard}. Four shaping variables prescribe the planform of the wing, by changing the root chord, the tip chord, the semispan, and the wing sweep. Each of these variables ranges from -5 to 5 inches, and additively scales the baseline wing shape. The baseline root chord is 22 inches, and so may vary during design from 17 to 27 inches. Similarly, the baseline tip chord is 14.5 inches, the baseline semispan is 30 inches, and the baseline wing sweep is 32 inches. In addition, we consider five structural sizing variables governing the plate thickness down the span of the wing, as shown in the figure.  The spanwise thickness distribution is governed by a piecewise linear interpolation across the 5 control points.  The baseline plate thickness is 1 inch, and may vary between 0.25 and 1.75 inches.

The aeroelastic constraints are detailed next. The function $r_1(\theta,\delta)$ is a dynamic aeroelastic constraint, and entails computing the flutter instability point, namely the flutter dynamic pressure, via a matched point $p-k$ scheme \cite{vanZyl}. Hence, the constraint $r_1(\theta,\delta)\leq 0$ is satisfied when the flutter dynamic pressure exceeds a limiting value of 2 psi. The function $r_2(\theta,\delta)$ is a static aeroelastic constraint: using the wing angle of attack, the flexible wing is trimmed to a lift coefficient of 0.5 at a dynamic pressure of 1.5 psi. The elastic von-Mises stresses are then computed within the deformed wing, and these values are finally aggregated into a single scalar output via the Kreisselmeier-Steinhauser method \cite{Kreisselmeier}. Hence, the constraint $r_2(\theta,\delta)<0$ is satisfied when the stress aggregation function is less than the yield stress, i.e., all the finite elements are inside their failure envelope.  The drag coefficient needed for the objective function is computed from the same trimmed state.

The parameter $\delta\in  \mathbb{R}^{6}$ combines uncertain parameters and changing operating conditions. In particular, $\delta$ includes the Mach number of the static and dynamic aeroelastic physics, the mass-proportional Rayleigh damping coefficient, the stiffness-proportional Rayleigh damping coefficient, the kinematic viscosity of the flow, and the target lift coefficient for the trim state. These parameters vary over the hyper-rectangular set $\Delta=[0.4,0.9]\times[0,25]\times[0,0.001]\times[\num{0.8e-5},\num{1.2e-5}]\times[0.4,0.6]$ respectively. The scenarios used for training and testing are obtained by sampling from a uniform distribution supported in $\Delta$. 

The underlying structural wing model is idealized as a flat plate shell cantilevered along its root, and immersed in subsonic flow. The structure of the wing is modeled with a linear shell finite element model, and both the steady and unsteady aerodynamics are modeled with the linear doublet lattice method \cite{Blair}. A finite plate spline \cite{Appa} is used to pass wing deformations from the structure to the aerodynamics, as well as loads from the aerodynamics to the structure.  Given that only linear compressible aerodynamics are utilized here, this solver will become less accurate as the Mach number approaches unity, and aerodynamic nonlinearities become more prominent.  Nonlinear flow solvers, i.e., computational fluid dynamics, could be used for higher accuracy, but at a much greater computational cost \cite{Stanford}.

A single-point wing design is presented first. This design, denoted as $\theta_1^\star$, is obtained from (\ref{wc}) for ${\mathcal D}=\{\delta_{\text{nom}}\}$, where $\delta_{\text{nom}}$ is the geometric center of $\Delta$. Single-point designs might not exhibit sufficient robustness to uncertainty regardless of the choice of $\delta_{\text{nom}}$. Choosing $\delta_{\text{nom}}$ as one of the worst-case combination of uncertainties, a difficult point to identify given its non-trivial dependency on $\theta$ through the requirement functions, might not yield the intended result, i.e., we may obtain a conservative design having an overly large objective value or an insufficiently robust design. As such, multi-point approaches are preferable.  

\begin{figure}[htb]
 \centering
 \includegraphics[trim={4cm 6cm 4cm 4cm},clip,width=5.5in]{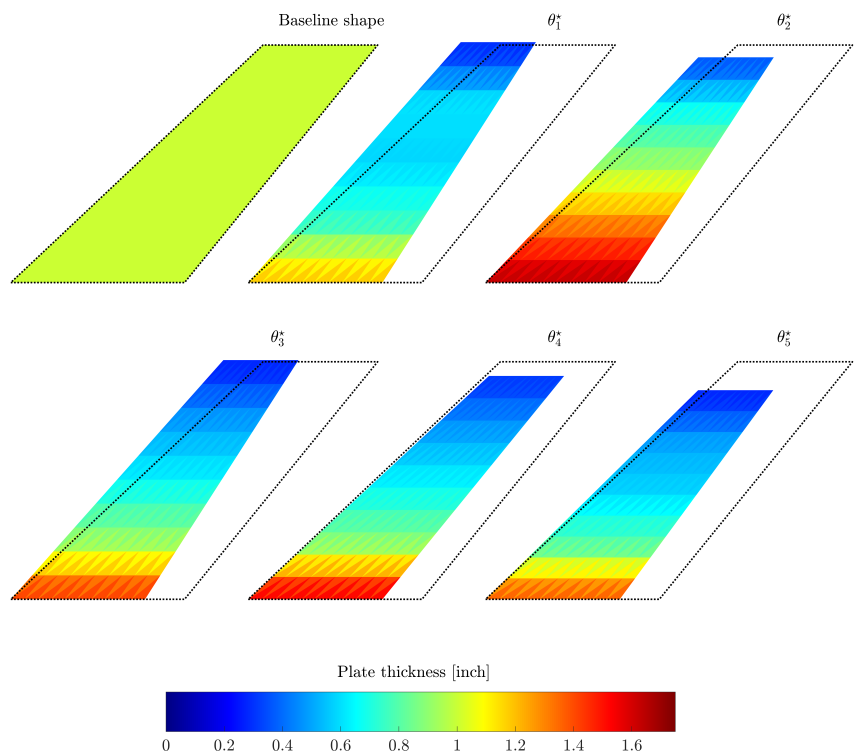}
 \caption{Optimal wing designs. The dashed black line indicates the baseline planform shape used for reference. }
 \label{wings}
\end{figure}

Two multi-point designs based on a training dataset ${\mathcal D}$ with $n=50$ scenarios and $m=1$ are presented next. A relatively small number of training scenarios is chosen due to the high computational cost of performing an aeroelastic analysis within the optimization loop. In particular, design $\theta_2^\star$ was computed using (\ref{chancon_sce}) with $J=0$, $\gamma_1=\gamma_2=1$. Hence, this design maximizes the probability of success for the nominal scenarios.  Design $\theta_3^\star$ was computed using (\ref{chancon_sce4}) for $\alpha=0$. Hence, this design minimizes the sample mean of the response for the nominal scenarios while ensuring that they satisfy the requirements. Figure \ref{wings} shows these and other wing designs. 

Prominent features of the resulting designs are discussed next. In an effort to minimize the structural mass, each design decreases the root and tip chord length, and also decreases the structural thickness at the tip.  A tapered thickness profile from root to tip is driven by a reduction in static stresses down the span.  Each optimal wing design is unswept relative to the baseline design, shown to the left of Figure \ref{wings}, a change that helps to satisfy the flutter constraint without increasing the structural mass.  The multi-point designs $\theta_2^\star$ and $\theta_3^\star$ exhibit  a large increase in the structural root thickness, which in turn lower the failure probability relative to $\theta_1^\star$ while increasing the objective value.

A reliability analysis of the above wing designs is presented next. This analysis, whose results are shown in Table \ref{tab0}, is carried out by using a Monte Carlo campaign with $n'=\num{1e4}$ sample points. Note that the 2nd, 3rd and 4th columns correspond to the training dataset ${\mathcal D}$ whereas the last six columns correspond to the testing dataset ${\mathcal D}'$. As expected, the large failure probabilities of $\theta_1^\star$ make it unsuitable.  That is not the case for $\theta_2^\star$ which meets the requirement for all scenarios in both the training and testing datasets. Note that a failure probability of zero based on $n'=\num{1e4}$ testing scenarios was attained by only using 50 training scenarios.  This outcome, however, is fortuitous. Recall that $\theta_2^\star$ ignores the response function but instead focuses on satisfying the requirements for all scenarios, thereby leading to a comparatively large performance value, ${\mathbb E}_\delta[h]$. This is not the case for $\theta_3^\star$, which lowers such value to 3.25 while satisfying the requirements for the training set. Note however that the reliability analysis of $\theta_3^\star$ reveals failure probabilities greater than zero. This outcome is caused by such a design overfitting the training dataset. 

\begin{table}[h!]
\scriptsize
\caption{Reliability analysis and supporting information: $n$ is the number of scenarios used for training, $m$ is the number of perturbations per scenario, $\sigma$ is the number of outliers, ${\mathbb E}_{_\delta}[h]$ is the empirical mean of the response function, $P_\delta[\mathcal{F}_i]$ is the failure probability for the $i$th requirement, $P_\delta[\mathcal{F}]$ is the total failure probability, and $\ell_i={\mathbb E}[r_i(\theta,\delta)|r_i(\theta,\delta)>0]$ is a loss measure. The last six columns use a testing dataset ${\mathcal D}'$ having $n'=\num{1e4}$ samples.}
\centering
 \label{tab0}
 \begin{tabular}{|c| c c c c | c c c c c c|}
 \hline 
               &      $n$   &      $m$         &         $\sigma$            &   ${\mathbb E}_{_\delta}[h]$   & ${\mathbb E}_\delta[h]$     & $P_\delta[\mathcal{F}_1]$              & $P_\delta[\mathcal{F}_2]$       & $P_\delta[\mathcal{F}]$    & $\ell_1$    & $\ell_2$\\ [0.5ex]
\hline
$\theta_1^\star$    &  1    &  1       &  0                             &3.08                                &    3.12                        &  \num{3.75e-1}                   &    \num{4.81e-1}                             &   \num{6.73e-1}  & 0.072 & 0.127 \\[0.5ex] 
$\theta_2^\star$    &  50  &  1       &  0                             &3.80                                &    3.81                        &  \num{0}                            &    \num{0}                                        &   \num{0}  & 0 & 0  \\[0.5ex] 
$\theta_3^\star$    &  50  &  1       &  0                            &3.25                                &    3.25                        &  \num{9.20e-3}                   &    \num{6.90e-3}                              &   \num{1.60e-2}  & 0.011& 0.012  \\[0.5ex] 
$\theta_4^\star$    &  50  &  $[1,2]$  &  0                        &3.31                               &    3.31                       &  \num{3.70e-3}                             & \num{1.00e-4}                          &  \num{3.80e-3}    &  0.088&  0.001  \\[0.5ex] 
$\theta_5^\star$    &  67  &  $1$  &  0                        &3.33                               &    3.28                       &  \num{0}                             & \num{0}                          &  \num{0}    &  0&  0  \\[0.5ex] 

\hline
\end{tabular}
\end{table}

The results of a robustness analysis, not shown here due to space limitations, are consistent with those of the reliability analysis. The wing design exhibiting the desired trade-off between performance and reliability/robustness should be chosen among the available alternatives.

Next we pursue a multi-point-robust design by using the adversarial perturbations in Section \ref{advper} for a small subset of the nominal scenarios. The corresponding sequence of adversarial perturbations will be denoted as ${\mathcal D}_p(\mathcal{Q}, \hat{k}, \hat{\theta})$, where $\mathcal{Q}=\{q^{(i)}\}_{i=1}^{n}$ and $q^{(i)}\in\{0,1\}$ is the number of times the $i$th scenario is perturbed. Formulation (\ref{chancon_sce4}) with $\alpha=0$ was used to synthesize the multi-point-robust wing $\theta_4^\star$ using the training dataset ${\mathcal D}_p(\mathcal{Q}, \hat{k}, \theta_3^\star)$. Whereas 46 elements of $\mathcal{Q}$ take the value of zero, four elements take the value of one. Therefore, only four scenarios are perturbed from their nominal value. These four scenarios attain the greatest $r_{\hat{k}}(\theta_3^\star,\delta^{(i)})$ negative values. Therefore,  $\theta_4^\star$ seeks to improve the robustness of $\theta_3^\star$ without significantly increasing the computational cost required for its calculation. Figures of merit corresponding to wing design $\theta_4^\star$, whose geometry is shown in Figure \ref{wings}, are listed in Table \ref{tab0}. This wing not only satisfies the requirements for the training dataset but also exhibits better robustness properties than $\theta_3^*$. In particular, $\theta_4^\star$ reduces the failure probability by a factor of 4.88 in exchange for an objective value increase of less than $2\%$. This robustness improvement is attained by only using four additional training scenarios. The loss measure of the flutter instability, as measured by $\ell_1$, increases by a factor of eight, whereas the loss measure for the stress, as measured by $\ell_2$, decreases by a factor of twelve.

New we use the sequential design strategy outlined in Section \ref{seque}. The process used to compute $\theta_4^\star$ constitutes the first iteration. In the second iteration we use formulation (\ref{chancon_sce2}) with $m=1$, $\gamma_1=\gamma_2=1$ for a dataset with $n=67$ scenarios. The 17 scenarios added to the original sequence of 50 scenarios fall onto ${\mathcal F}(\theta_4^\star)$ while attaining the smallest $\min_k r_k(\theta_4^\star,\delta^{(i)})$ values. The resulting design, denoted as $\theta_5^\star$, drives the failure probability to zero while attaining a much smaller objective value than $\theta_2^\star$ (Table \ref{tab0}), the other design attaining failure probabilities of zero.  Figure \ref{wings} shows that the most robust designs, $\theta_2^\star$ and $\theta_5^\star$, differ significantly. 

\begin{figure}[htb]
 \centering
 \includegraphics[trim={0.4cm 0cm 2cm 0cm},clip,width=4in]{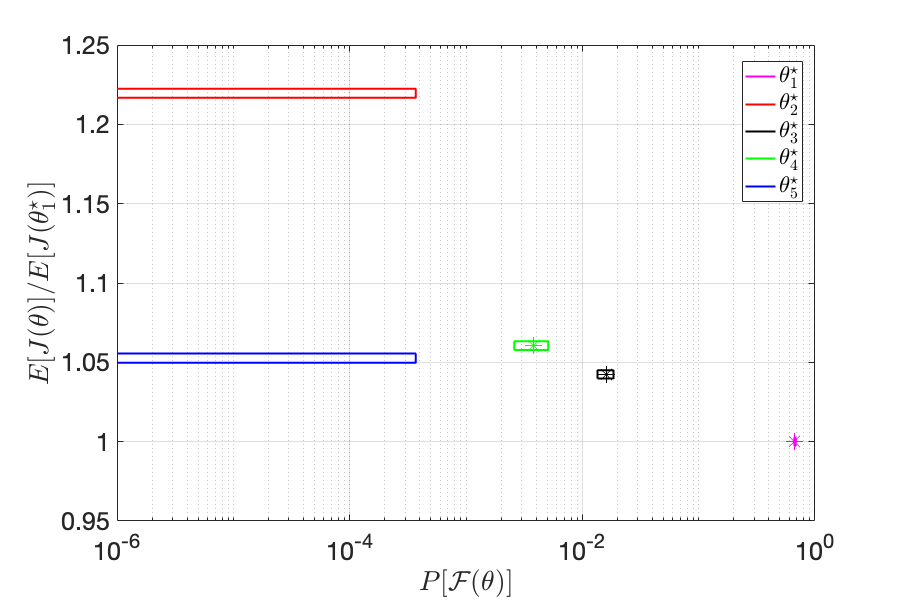}
 \caption{95 confidence intervals of the probability of failure and the expected response based on a testing dataset with $n'=\num{1e4}$ samples. The empirical estimates corresponding to non-zero failure probabilities are shown as asterisks.}
 \label{compare}
\end{figure}

Figure \ref{compare} shows the probability of failure against the normalized expected response for the 5 wing designs. This figure not only shows the empirical estimates corresponding to the testing dataset but also the corresponding 95\% confidence intervals.  The volume of these intervals, which lead to the rectangles shown, approaches zero as $n'$ approaches infinity. Note that both $\theta_2^\star$ and $\theta_5^\star$ drive the failure probability to zero but the latter does it by only increasing the objective value 6\%.  This illustrates the potential drawbacks of designs focusing on satisfying the requirements only. Monte Carlo campaigns with a greater number of samples are required to reduce the width of the confidence intervals, thereby further discriminating $\theta_2^\star$ from $\theta_5^\star$. The high computational cost of such a simulation, however, might render this practice infeasible.

The optimization programs above were solved using the active-set algorithm from Matlab, i.e., a sequential quadratic programming method, without analytical derivatives. The calculation of an optimal design required a few thousand iterations to converge when it started from the baseline wing design. This task took a couple of days of CPU time in a standard personal computer. 

\section{Concluding Remarks} 
This paper proposes a scenario optimization framework to robust design in the presence of error and uncertainty in the data. The system's performance, as measured by the value taken by the objective function, and the system's robustness, as measured by the ability to satisfy the requirements for perturbations in the data, are traded off by using two types of relaxations. Specifically, the feasible set is expanded by not only eliminating outliers from the dataset but also replacing constraints for the worst-case perturbation with chance constraints. Furthermore, we study the effects that loss measures commonly used to select outliers have on the resulting design, and propose formulations that do not depend on any of them. The optimization programs proposed, some of which have a number of decision variables that do not increase with the number of scenarios, can be solved using standard gradient-based algorithms while being applicable to continuous but otherwise arbitrary requirement functions. This setting is amenable to many problems in science and engineering for which the objective function and constraints must be evaluated by simulation.  Furthermore, we propose using targeted, adversarial perturbations and sequential data-driven design to mitigate the computational cost of performing optimization under uncertainty. The strategies proposed can be naturally integrated to the Monte Carlo campaigns commonly used to evaluate a system's robustness by bridging the verification phase with the design phase. 



\section*{Appendix (CDF approximations)}
Approximations to the CDF of a random variable and its inverse are presented next. Consider the non-decreasing sequence ${\mathcal Z}=\{z_i\}_{i=1}^{n}$ that results from evaluating the function $z(\theta,\delta)$ at $\{\delta^{(i)}\}_{i=1}^n$ for a fixed $\theta$ and sorting the resulting values so $z_i=z(\delta^{(j)})$ for some $j$ and $z_i < z_{i+1}$.  A continuous, piecewise linear approximation to the CDF of $z(\delta)$ based on ${\mathcal Z}=\{z(\delta^{(i)})  \}_{i=1}^n$ is
\begin{align}\label{Flin}
F_{{\mathcal Z}(\theta)}(z)\triangleq\begin{cases}
0 & \text{if } z\leq z_1,\\
\frac{1}{n-1}\left(i-1 + \frac{z-z_i}{z_{i+1}-z_i}\right) & \text{if } z_i < z \leq z_{i+1}, \\
1 & \text{otherwise.}
\end{cases}
\end{align}
The inverse of (\ref{Flin}) is
\begin{align}\label{Finvlin}
F^{-1}_{{\mathcal Z}(\theta)}(\alpha)=\begin{cases}
z_1 & \text{if } \alpha=0,\\
z_{i}+(z_{i+1}-z_i)\left((n-1) \alpha-i+1\right) & \text{if } 0<\alpha<1, \\
z_n & \text{otherwise,}
\end{cases}
\end{align}
where $\alpha\in[0,1]$ and $i=\argmin_{j=1,\dots n} \{(n-1)\alpha-j+1: j-1\leq \alpha(n-1)\}$. The approximations (\ref{Flin}) and (\ref{Finvlin}) are differentiable in $\theta$ when $z(\theta,\delta)$ is $C^1$ in $\theta$. This property makes standard gradient-based algorithms applicable to the above optimization programs. When ${\mathcal Z}$ contains repeated values, (\ref{Flin}) and (\ref{Finvlin}) can be used after breaking the ties with small perturbations.


\section*{Acknowledgements}
This work was supported by the NASA Human Research Program (HRP) for radiation protection.



 \bibliographystyle{elsarticle-num} 
 
 {\footnotesize \bibliography{nref}}





\end{document}